\documentclass[authoryear,11pt,preprint]{elsarticle}
\usepackage{color}
\usepackage{lineno}
\usepackage{pdflscape}
\usepackage{amsmath}
\usepackage{booktabs}
\usepackage[hyphens]{url}
\usepackage{hyperref}
\usepackage{multirow}
\usepackage{addlines}
\usepackage{adjustbox}
\usepackage{natbib} %When changing reference style
\usepackage[english]{babel}
\usepackage[dvipsnames]{xcolor}
\usepackage{soul}
\usepackage{adjustbox}
\usepackage{verbatim}
\usepackage[margin=2.5cm]{geometry}% by courtesy of Mico
\usepackage{mathtools}
\usepackage{tikz}
\usepackage{pgfplots}
\usepackage{enumitem}
\usepackage{calc}

\definecolor{darkblue}{rgb}{0,0.2,0.4}
\usepackage{graphicx}
\usepackage{subcaption}
\usepackage{amssymb}
\usepackage{amsthm}
\usepackage{etoolbox}
\usepackage{algorithm}
\usepackage{algpseudocode}

\newtheorem{cond}{Condition}
\newtheorem{proper}{Property}
\newtheorem{theorem}{Theorem}
\newtheorem{lemma}{Lemma}

\usepackage{float} % For defining new floating environments
\usepackage{algpseudocode} % For algorithmic structures
\usepackage{algorithm} % For algorithm floating environment styling

\newfloat{subroutine}{htbp}{loa}
\floatname{subroutine}{Subroutine}

\makeatletter
\newcommand*{\algrule}[1][\algorithmicindent]{%
  \hspace*{.2em}% <------------- This is where the rule starts from
  \vrule %height .75\baselineskip depth .25\baselineskip
  \hspace*{\dimexpr#1-.2em-.4pt}%
}

\newcommand{\StatePar}[1]{%
  \State\parbox[t]{\dimexpr\linewidth-\ALG@thistlm}{\strut #1\strut}%
}
\renewcommand{\ALG@beginalgorithmic}{\offinterlineskip}% Remove all interline skips

\newcount\ALG@printindent@tempcnta
\def\ALG@printindent{%
  \ifnum \theALG@nested > 0% is there anything to print
    \ifx\ALG@text\ALG@x@notext% is this an end group without any text?
    \else
      \unskip
      \ALG@printindent@tempcnta=1
      \loop
        \algrule[\csname ALG@ind@\the\ALG@printindent@tempcnta\endcsname]%
        \advance \ALG@printindent@tempcnta 1
        \ifnum \ALG@printindent@tempcnta<\numexpr\theALG@nested+1\relax
      \repeat
        \fi
    \fi
}
\patchcmd{\ALG@doentity}{\noindent\hskip\ALG@tlm}{\ALG@printindent}{}{\errmessage{failed to patch}}
\makeatother

\algrenewcommand\algorithmicend{\strut\textbf{end}}
\algrenewcommand\algorithmicdo{\strut\textbf{do}}
\algrenewcommand\algorithmicwhile{\strut\textbf{while}}
\algrenewcommand\algorithmicfor{\strut\textbf{for}}
\algrenewcommand\algorithmicforall{\strut\textbf{for all}}
\algrenewcommand\algorithmicloop{\strut\textbf{loop}}
\algrenewcommand\algorithmicrepeat{\strut\textbf{repeat}}
\algrenewcommand\algorithmicuntil{\strut\textbf{until}}
\algrenewcommand\algorithmicprocedure{\strut\textbf{procedure}}
\algrenewcommand\algorithmicfunction{\strut\textbf{function}}
\algrenewcommand\algorithmicif{\strut\textbf{if}}
\algrenewcommand\algorithmicthen{\strut\textbf{then}}
\algrenewcommand\algorithmicelse{\strut\textbf{else}}
\algrenewcommand\algorithmicrequire{\strut\textbf{Input:}}
\algrenewcommand\algorithmicensure{\strut\textbf{Output:}}
\let\oldState\State
\renewcommand{\State}{\oldState\strut}
\algnewcommand{\IIf}[1]{\State\algorithmicif\ #1\ \algorithmicthen}
\algnewcommand{\EndIIf}{\unskip\ \algorithmicend\ \algorithmicif}

\usepackage{nomencl}
\makenomenclature
\usepackage{framed} % Framing content
\usepackage{multicol} % Multiple columns environment

\usepackage[utf8]{inputenc}
\usepackage{newunicodechar}

\DeclareUnicodeCharacter{0229}{\c{e}}
\DeclareUnicodeCharacter{2009}{ °}
\usepackage[T1]{fontenc}
\usepackage{textcomp}
\usepackage{lmodern}
\usepackage{ltxcmds}
\usepackage{gensymb}
\usepackage{siunitx}
\DeclareSIUnit{\euro}{\texteuro}
\makeatother
\definecolor{royalpurple}{rgb}{0.58, 0.44, 0.86}

\allowdisplaybreaks

\DeclareMathOperator*{\argmax}{argmax}
\DeclareMathOperator*{\argmin}{argmin}

\newcommand\norm[1]{\left\lVert#1\right\rVert}

\usepackage{setspace}
\hypersetup{colorlinks = true,linkcolor = blue,anchorcolor =red,citecolor = blue,filecolor = red,urlcolor = red, pdfauthor=author}

\journal{}

\pgfplotsset{compat=1.17} % for plotting
\begin{document}

\begin{frontmatter}

\title{Integrated investment, retrofit and abandonment energy system planning with multi-timescale uncertainty using stabilised adaptive Benders decomposition}


\author[a,b]{Hongyu Zhang\corref{cor1}}
\ead{hongyu.zhang@soton.ac.uk}
\author[c]{Ignacio E. Grossmann}
\ead{grossmann@cmu.edu}
\author[d]{Ken McKinnon}
\ead{K.McKinnon@ed.ac.uk}
\author[e]{Brage Rugstad Knudsen}
\ead{brage.knudsen@sintef.no}
\author[d]{Rodrigo Garcia Nava}
\ead{rodrigo.garciana@gmail.com}
\author[a]{Asgeir Tomasgard}
\ead{asgeir.tomasgard@ntnu.no}

\cortext[cor1]{Corresponding author}

\address[a]{Department of Industrial Economics and Technology Management, Norwegian University of Science and Technology, Høgskoleringen 1, 7491, Trondheim, Norway}
\address[b]{School of Mathematical Sciences, University of Southampton, Building 54, Highfield Campus, Southampton, SO14 3ZH, United Kingdom}
\address[c]{Department of Chemical Engineering, Carnegie Mellon University, 5000 Forbes Avenue, Pittsburgh, PA 15213, USA}
\address[d]{School of Mathematics, University of Edinburgh, James Clerk Maxwell Building, Peter Guthrie Tait Road, Edinburgh, EH9 3FD, United Kingdom}
\address[e]{SINTEF Energy Research, Kolbjørn Hejes vei 1B, 7491, Trondheim, Norway}

\begin{abstract}

We propose the REORIENT (REnewable resOuRce Investment for the ENergy Transition) model for energy systems planning with the following novelties: (1) integrating capacity expansion, retrofit and abandonment planning, and (2) using multi-horizon stochastic mixed-integer linear programming with multi-timescale uncertainty. We apply the model to the European energy system considering: (a) investment in new hydrogen infrastructures, (b) capacity expansion of the European power system, (c) retrofitting oil and gas infrastructures in the North Sea region for hydrogen production and distribution, and abandoning existing infrastructures, and (d) long-term uncertainty in oil and gas prices and short-term uncertainty in time series parameters. We utilise the structure of multi-horizon stochastic programming and propose a stabilised adaptive Benders decomposition to solve the model efficiently. We first conduct a sensitivity analysis on retrofitting costs of oil and gas infrastructures. We then compare the REORIENT model with a conventional investment planning model regarding costs and investment decisions. Finally, the computational performance of the algorithm is presented. The results show that: (1) when the retrofitting cost is below 20\% of the cost of building new ones, retrofitting is economical for most of the existing pipelines, (2) platform clusters keep producing oil due to the massive profit, and the clusters are abandoned in the last investment stage, (3) compared with a traditional investment planning model, the REORIENT model yields 24\% lower investment cost in the North Sea region, and (4) the enhanced Benders algorithm is up to 6.8 times faster than the level method stabilised adaptive Benders.

\end{abstract}

\begin{keyword}
OR in energy \sep Stochastic programming \sep Multi-horizon stochastic programming \sep Large-scale mixed-integer linear programming \sep Retrofit of energy systems

% MSC codes here, in the form: \MSC code \sep code
% or \MSC[2008] code \sep code (2000 is the default)

\end{keyword}
%Highlights
% Presemt a new vision on creating vir
\end{frontmatter}

% \linenumbers

\section{Introduction}
\label{sec:introduction}
Accelerating energy transition in all sectors is vital to achieve a carbon-neutral economy by 2050 \citep{REPowerEU, green_deal}. The committed emissions from existing energy infrastructure jeopardise the 1.5 °C target \citep{Tong2019CommittedTarget}. It may be more beneficial to retrofit existing energy infrastructure than to abandon it. Abandoning existing energy infrastructure, such as oil and gas infrastructure, may have a substantial cost \citep{Bakker2019AnCampaigns}. Also, the oil and gas industry involves multi-billion-dollar investments and profits. Therefore, there is motivation to retrofit existing oil and gas infrastructure for clean energy production and transportation to (1) help the oil and gas industry transition to a clean energy producer, and (2) accelerate the energy transition by financing it using the gain from the avoided abandonment cost. Retrofitting existing oil and gas infrastructure for hydrogen production and transportation is drawing more attention due to the increasing demand for hydrogen. Most offshore pipelines can be used for hydrogen transport in Europe \citep{ValentinGentile2021Re-StreamEurope}. The European hydrogen infrastructure could grow to a pan-European network by 2040, largely based on repurposed existing natural gas infrastructure \citep{RikvanRossum2022EuropeanBackbone}. Retrofitting an existing offshore platform to a green hydrogen production platform is under exploration \citep{NeptuneEnergy2023PosHYdonSea}. We note that retrofitting may become an important pillar for accelerating the energy transition.  Therefore, in this paper, we analyse cost drivers that have the possibility to trigger widespread retrofit of offshore oil and gas infrastructure and achieve cleaner energy generation and decarbonisation. We use a high-fidelity, detailed spatial-temporal stochastic programming model to analyse these drivers for a large region with a substantial role in the energy supply to the surrounding countries.

Energy system infrastructure planning is crucial during the energy transition towards zero emission by $2050$. Optimisation models are widely used for the investment \citep{zhang2022_OEH,Cho2022RecentOptimization, Munoz2015} and operational \citep{Schulze2016,Philpott2000Hydro-electricDemand} planning of energy systems. Traditionally, capacity expansion, retrofitting \citep{store2018} and abandonment are planned separately using different models. However, as sector coupling and Power-to-X become more important, as well as the possibility of retrofitting existing infrastructure for renewable energy production and distribution, the ability to analyse investments, retrofit and abandonment planning in a single integrated model becomes more important; including all degrees of freedom together minimises the risk of suboptimality. However, such integrated models have been less explored than their counterparts that treat investment, retrofit and abandonment independently. 

Managing uncertainty is crucial in a long-term planning model. Long-term energy system planning problems are normally prone to uncertainty on strategic and operational time horizons. Strategic uncertainty includes oil and gas prices, CO$_2$ budget, and CO$_2$ tax. Short-term uncertainty normally includes the availability of non-dispatchable renewable technologies. When the short-term uncertainty does not affect the future strategic decisions, the problem has a Multi-Horizon Stochastic Programming (MHSP) structure \citep{Kaut2014}, and this leads to a more compact and easier to solve model than when this is not the case. Most previous studies on energy system planning consider only short-term uncertainty \citep{Backe2022EMPIRE:Analyses}. In this paper, the proposed model utilises MHSP and includes uncertainty from both time horizons.

Large-scale MHSP remains difficult, however, and this needs to be addressed. The block separable structure of MHSP allows the decomposition of a problem with Benders-type algorithms. Also, the structure of the subproblems enables adaptive oracles \citep{Mazzi2020}. A stabilised adaptive Benders decomposition algorithm was proposed in \cite{zhang2022_benders} and demonstrated on a power system investment planning problem with up to 1 billion variables and 4.5 billion constraints. The algorithm showed a significant reduction in computational time. However, \cite{zhang2022_benders} dealt with a large-scale linear programming problem. In this paper, we consider a problem with binary variables in the investment planning part in order to capture the economic scale and model retrofit and abandonment decisions. Because the binary variables only exist in the investment planning part, which is the reduced master problem in the Benders-type algorithm, we can still apply stabilised adaptive Benders directly. However, the algorithm may slow down due to the combinatorial part of the problem. This is because the stabilisation problem in \cite{zhang2022_benders} are quadratic programs, and when binary variables are added, the problems become mixed integer quadratic programs. Although commercial solvers have improved, allowing mixed integer quadratic programming to be solved more efficiently, they are unable to deal with the large problems this paper considers.

To fill the research gaps mentioned above, we first propose the REORIENT (REnewable resOuRce Investment for the ENergy Transition) model, which is a multi-horizon stochastic Mixed-Integer Linear Programming (MILP) model with multi-timescale uncertainty for the investment, retrofit and abandonment planning for energy systems.  We consider retrofitting existing platforms for offshore green hydrogen production and pipelines for green and blue hydrogen distribution. We then demonstrate the REORIENT model on a European energy system planning problem. A stabilised adaptive Benders decomposition is proposed: this extends and improves the method in \cite{zhang2022_benders}, allowing it to solve problems with binary variables efficiently.

The contributions of the paper are the following: (1) we integrate investment planning, retrofitting, and abandonment in a single multi-horizon stochastic MILP model, which allows us to analyse the economic trade off among different kinds of planning decisions rather than only capacity expansion decisions in a traditional investment planning model, (2) we propose a centred point stabilised Benders decomposition with inexact cuts to solve the resulting large-scale model, (3) we extend the centred point stabilisation, which was proposed for linear programs, to solve MILP and provide convergence proof, and (4) we demonstrate the model on a large-scale planning problem for the European energy system to analyse the planning decisions and costs and provide global insights.

The outline of the paper is as follows: Section \ref{sec:literature_review} introduces the background knowledge regarding capacity expansion planning, retrofit planning and abandonment planning, stochastic programming, MHSP, and Benders decomposition. Section \ref{sec:problem_description_modelling_assumptions} provides the problem description, modelling strategies and assumptions. Section \ref{sec:Benders decomposition in MHSP} presents the proposed stabilised adaptive Benders decomposition. Section \ref{sec:model} presents the REORIENT model. Section \ref{sec:results} reports the computational results and numerical analysis. Section \ref{sec:discussion} discusses the implications of the method and results and summarises the limitations of the research. Section \ref{sec:conclusions} concludes the paper and suggests further research.

\section{Literature review}
\label{sec:literature_review}
In the following, we present a  brief overview of relevant literature on capacity expansion planning, abandonment planning, retrofit planning, MHSP, and Benders decomposition.

\subsection{Capacity expansion planning}
\label{sec:capacity expansion planning}
Capacity expansion planning problems normally consider an existing system with historical capacity or a new system and make investment planning to fulfil the demand under, among others, environmental restrictions. Capacity expansion problems are formulated either in deterministic models \citep{Lara2018DeterministicAlgorithm} or stochastic programming models \citep{Backe2022EMPIRE:Analyses, Conejo2016}. \cite{Backe2022EMPIRE:Analyses} utilised multi-horizon formulation but did not include long-term uncertainty. Sometimes MILP is used \citep{lara2020} to capture the discrete nature of some investment decisions. To gain enough environmental and economic insights from such models, sometimes large-scale problems need to be solved, such as in \cite{Li2022} and \cite{zhang2022_OEH}. \cite{Munoz2016} proposed a new bounding scheme and combined it with Benders decomposition to solve a large investment planning problem that was formulated as MILP. In addition to planning for power, natural gas and heat systems separately, planning for integrated multi-carrier systems has also been studied. Energy hubs that convert, process and store multiple energy carriers in an investment planning problem were studied in \cite{zhang2022_benders}. Also, stochastic programming has been used in natural gas systems \citep{Fodstad2016}, offshore oil and gas infrastructure planning \citep{Gupta2014}, and hydrogen network \citep{Galan2019}. There is much literature on capacity expansion problems: we refer the readers to \cite{Krishnan2016Co-optimizationApproaches} for a more comprehensive review.

\subsection{Retrofit planning}
\label{sec:retrofit planning}
In grassroots design, the decisions on what processes to use are made first to be followed by equipment decisions, but for proper analysis the retrofit design also requires models that evaluate existing equipment. A comparison of grassroots and retrofit design has been presented in \cite{grossmann1987}. The combinatorial nature of the retrofit planning problems makes their models much more complex. There are several reasons to retrofit, including (1) processing a new feedstock, (2) improving the economics by the use of less energy per unit of production, (3) making a new product, and (4) increasing the conversion of the feedstocks. In this paper, we consider retrofitting to make a new product. A debottlenecking strategy was proposed for retrofit problems \citep{Harsh1989AProblem}. A systematic procedure for the retrofit of heat exchanger networks was presented in \cite{yee1991}. Retrofitting of heat exchangers has been extensively studied in the past decades \citep{Pan2013ADetails,Wang2012ApplicationNetwork}. A high-level optimisation model for the retrofit planning of process networks over several time periods was presented in \cite{Jackson2002High-levelNetworks}. The proposed strategy consists of a high-level analysis of the entire network and a detailed low-level analysis of the resulting specific process flowsheet. The problem is formulated using a multi-period generalised disjunctive programming model, which is reformulated using a mixed-integer linear program using the convex hull formulation. In this paper, because of the scope and size of the problem, we only consider the high-level modelling of the retrofit. 

\subsection{Abandonment planning}
\label{sec:plug and abandonment planning}
Abandonment planning includes the abandonment of power plants that exceed their lifetimes, and of mature oil and gas fields. In the oil and gas industry, existing literature focuses on the plug and abandonment campaign. This is because many wells are planned to be plugged and abandoned, and such activity will have a substantial cost \citep{Bakker2019AnCampaigns}. The plug and abandonment cost is estimated at £5-15 million per well, and thousands are expected to be abandoned just in offshore regions over the next decade. Plug and abandonment planning is formulated either as a profit maximisation problem \citep{Bakker2021MatureProgramming} or a cost minimisation problem \citep{Bakker2019AnCampaigns, Bakker2021VehiclePlanning}. Real options theory is also used for oil and gas field development \citep{Stein-ErikFleten2011RealTie-ins,store2018,Bakker2019AnCampaigns}. In addition to abandonment, investment and operational scheduling optimisation in the oil and gas sector can be found in \cite{Oliveira2013AConsiderations,Goel2004AReserves,Gupta2012AnInfrastructure,iyer1998}.

\begin{figure}[!htb]
    \centering
    \includegraphics[scale=0.8]{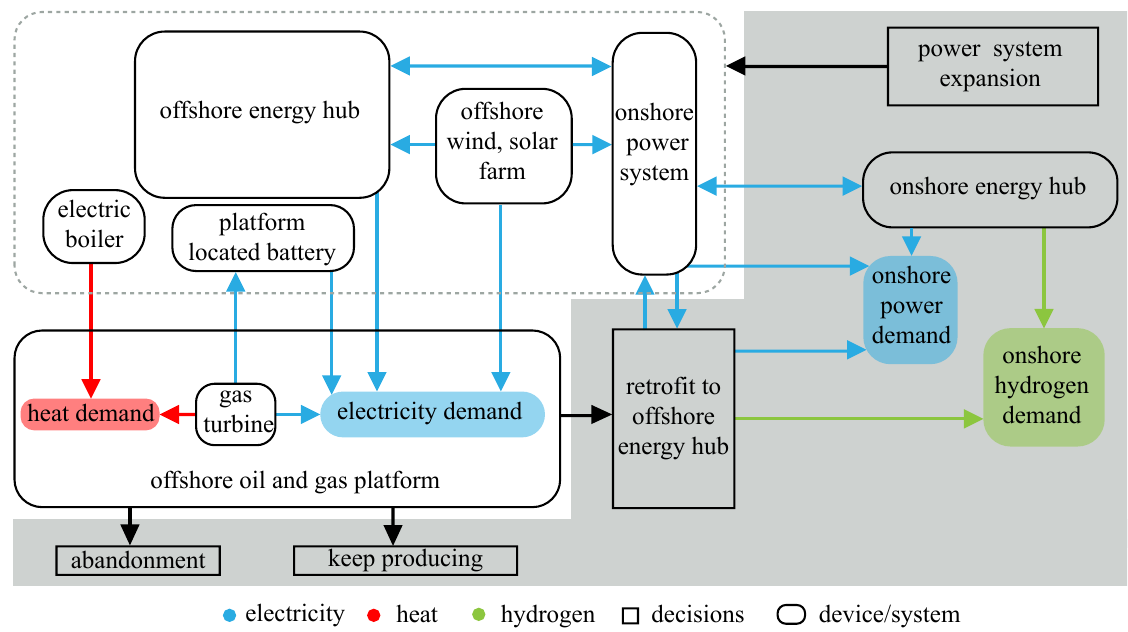}
    \caption{Integrated investment, retrofit and abandonment planning. The grey dotted box includes some technologies that can be invested in. The offshore oil and gas platform can be retrofitted or abandoned. Otherwise, it can keep producing. The Offshore Energy Hubs (OEHs) \citep{zhang2022_OEH} can generate electricity, and produce and store hydrogen.}
    \label{fig:retrofit}
\end{figure}

From the literature above, we find that the planning of investment, retrofit, and plug and abandonment are often treated separately. To overcome the limitations of the separated approach for energy systems planning where such decisions are deemed tightly coupled, we propose a modelling framework that integrates investment, plug and abandonment and retrofit. An illustrative example is presented in Figure \ref{fig:retrofit}. The parts marked with grey represent the new integrated planning compared with traditional investment planning in the literature. 

\subsection{MHSP}
\label{sec:multi_horizon_programming}
Investment planning of an energy system often faces uncertainty from two time horizons \citep{Kaut2014,lara2020}: (a) the uncertainty from the operational time horizon, such as the availability of generation from renewable energy sources. The operational uncertainty becomes even more crucial for a system with higher penetration of intermittent renewable energy, and (b) the uncertainty on the strategic time horizon, e.g., oil and gas price and demand. In traditional multi-stage stochastic programming, uncertainty from operational and strategic time horizons can lead to a large scenario tree, thus, an intractable planning model. The multi-horizon approach was proposed as an alternative formulation that reduces the model size significantly by reducing the interaction between short-term and long-term uncertainties \citep{Kaut2014}.

\begin{figure*}[!htb]
    \centering
    \includegraphics[scale=0.8]{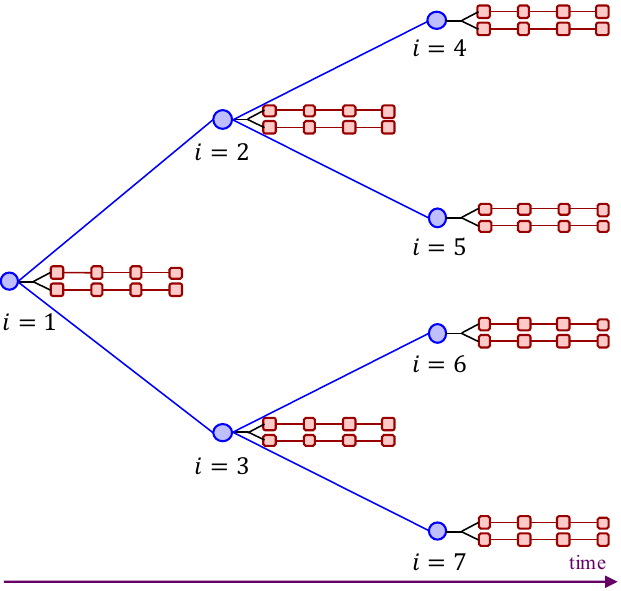}
    \caption{Illustration of MHSP with short-term and long-term uncertainty. (blue circles: strategic nodes, red squares: operational periods, $i:$ index of the strategic nodes)}
    \label{fig:MHSP}
\end{figure*}

One can have a much smaller model by disconnecting operational nodes between successive planning stages and embedding them into their respective strategic nodes. An illustration of MHSP with short-term and long-term uncertainty is shown in Figure \ref{fig:MHSP}. We call an operational problem embedded in a strategic node an operational node. The resulting model is called MHSP. MHSP is identical to multi-stage stochastic programming provided that the state of an operational problem at the end of a strategic period does not affect the operational state at the start of the following strategic period or affect any future strategic decisions. If either of these conditions is not met, then MHSP gives only an approximation, however for our problem it is a good approximation. MHSP has been applied in several energy system planning problems \citep{Skar2016,Backe2022EMPIRE:Analyses,Zhang2022OffshoreShelf,Durakovic2023PoweringPrices}. Furthermore, the bounds in MHSP have been studied in \cite{Maggioni2020BoundsPrograms}. 

\subsection{Benders decomposition}
Benders decomposition was proposed to solve problems with complicating variables \citep{Geoffrion1972}. Then generalised Benders for nonconvex problems was proposed \citep{Steeger2017DynamicProblem,Li2019AVariables}. Extensive research was conducted on accelerating Benders decomposition, such as by choosing and adding strong cuts \citep{Magnanti1981AcceleratingCriteria, Oliveira2014AcceleratingChain}, developing a multi-cut version of Benders decomposition \citep{You2011MulticutUncertainty}, stabilised Benders decomposition \citep{Zverovich2012,zhang2022_benders}, and using inexact oracles \citep{Mazzi2020,zhang2022_benders,vanAckooij2016InexactSupport}. In this paper, we adopt the multi-cut approach used in \cite{zhang2022_benders} and improve and extend the method to solve MILPs. 

A multi-stage stochastic capacity expansion program formulated as MHSP has a block separable structure \citep{Louveaux1986} and can be decomposed by Benders decomposition \citep{zhang2022_benders}.  For Benders decomposition, all the long-term nodes (blue circles in Figure  \ref{fig:MHSP}) are included in a single reduced master problem, and the blocks of short-term nodes (red squares in Figure \ref{fig:MHSP}) are in the subproblems. The REORIENT model described in this paper is an MHSP problem, and the investment planning part is tractable even though all the long-term nodes are included in a single reduced master problem. Therefore, we choose to use and improve the method proposed in \cite{zhang2022_benders} instead of using nested Benders decomposition or stochastic dual dynamic programming, whose benefit is to have multiple but smaller reduced master problems.

A potential problem of Benders decomposition is that all subproblems need to be solved exactly at every iteration to generate cuts and obtain an upper bound. \cite{Mazzi2020} proposed Benders decomposition with adaptive oracles to avoid solving all subproblems at every iteration by utilising the structure of subproblems. In \cite{Mazzi2020}, Benders decomposition with adaptive oracles was shown to be 31.9 times faster than standard Benders decomposition for a 1\% convergence tolerance, and it has been further improved in \cite{zhang2022_benders}. The subproblems in the REORIENT model are large-scale operational problems that are hard to solve, and the adaptive oracles are needed to make the problem tractable.

\cite{zhang2022_benders} pointed out that stabilising Benders decomposition is important when solving multi-region energy system planning problems. The level method stabilised adaptive Benders decomposition proposed in \cite{zhang2022_benders} is able to solve linear programs with 1 billion variables and 4.5 billion constraints very efficiently. Stabilisation is needed since the REORIENT model is a large-scale multi-region planning model. In this paper, we propose a new stabilisation method for solving the REORIENT model. This is because the level method stabilisation problem in \citep{zhang2022_benders} experiences a severe slowdown when there are binary variables. Therefore, we propose an alternative stabilisation called centred point stabilisation that speeds up the stabilisation step.

In addition to Benders decomposition, Lagrangean decomposition \citep{Escudero2017ScenarioOptimization}, column generation \citep{Singh2009}, and combined column generation and Benders decomposition \citep{Huang2022AUncertainties} have been proposed for capacity expansion problems. These approaches can solve problems with integer variables in the operational problem, however, this is not a feature that is needed in the REORIENT model, and they do not exploit the block separable structure of MHSP, which makes them less suitable alternatives than the method in \cite{zhang2022_benders}.

\section{Problem description, modelling strategies, modelling assumptions}
\label{sec:problem_description_modelling_assumptions}
In this section, we present the problem description and modelling strategies, including price modelling, scenario generation, temporal and geographical representations of the problem, and the modelling assumptions. 

The problem under consideration aims to choose (a) the optimal strategy for investment, abandonment and retrofit planning, and (b) operating scheduling for an energy system to achieve emission targets at minimum overall costs under short-term uncertainty, including renewable energy availability, hydropower production profile and load profile, and long-term uncertainty, including oil and gas prices. 

For the investment planning, we consider: (a) thermal generators (Coal-fired plant, OCGT, CCGT, Diesel, nuclear plants, co-firing biomass with 10\% lignite, lignite); (b) generators with Carbon Capture and Storage (CCS) (Coal-fired plant with CCS and advanced CCS, gas-fired plant with CCS and advanced CCS, co-firing biomass with 10\% lignite with CCS, lignite with CCS, lignite with advanced CCS); (c) renewable generators (onshore and offshore wind and solar, wave, biomass, run-of-the-river hydropower, geothermal and regulated hydropower); (d) electric storage (hydro pump storage and lithium); (e) onshore and offshore transmission lines; (f) onshore and offshore clean energy hubs (electrolyser, fuel cell, hydrogen storage); (g) onshore steam reforming plant with CCS (SMRCCS) and (h) offshore and onshore hydrogen pipelines. The capital costs and fixed operational costs coefficients are assumed to be known. 

For the retrofit planning, we consider: (a) retrofitting existing natural gas pipelines for hydrogen transport, and (b) retrofitting existing offshore platforms for clean OEHs. Finally, we consider the abandonment of mature fields.  The problem is to determine: (a) the capacities of technologies and retrofit and abandonment decisions, and (b) operational strategies that include scheduling of generators, storage and approximate power flow among regions to meet the energy demand with minimum overall expected investment and operational and environmental costs.

\subsection{Modelling strategies and assumptions}
In this section, we present the modelling strategies and assumptions we use in the long-term integrated planning problem.

\subsubsection{Long-term production profile modelling}
There are three phases during the lifetime of oil field production \citep{store2018}: ramp up, plateau and decline. We adopt the commonly used production model from \cite{wallace1987}, whose long-term production quantity of oil and gas is represented by,

\begin{equation}
    q_t=\begin{cases}
        q^P & t\leq t^P,\\
        q^Pe^{-m(t-t^p)} &  t> t^P,
        \end{cases}
\end{equation}
where $q^P$ (MMbbl) is the production rate during plateau period, $t^P$ (year) is the length of the plateau period, $m$ is the decline constant. We also calibrate this model using the average decline rate of the giant oil fields in the world \citep{hook2009}. Note that with this model, the amount of oil and gas in a reservoir at the start of each planning stage, and hence the total amount produced within a planning stage, are exogenous quantities, and so are not affected by the operations. This approximation maintains the MHSP structure of the model.

\subsubsection{Scenario generation}
For short-term uncertainty, we consider uncertain time series data: wind and solar capacity factors, hydropower generation profiles, load profiles, and platform production profiles. We consider operational problems with an hourly resolution. We divide the full year into four seasons and select representative time slices from these seasons. The length of the times slices can be different in different seasons. The spatial correlation is preserved due to sampling from historical data. To preserve the auto-correlation and correlation between time series data, the same hours are used for all the time series data within a scenario. The generated operational scenarios are used for all operational nodes.

For long-term uncertainty, we consider oil and gas price uncertainty. We first generate a large number of projections using the price process described below. We use a two-factor Short-Term Long-Term (STLT) model to capture a realistic price behaviour of oil and gas \citep{Schwartz2000Short-termPrices}. The STLT model has a stochastic equilibrium price and stochastic short-term deviations \citep{Bakker2021MatureProgramming}. In the STLT model, the logarithm of the spot price at time $t$ is, 

\begin{equation}
    \log p_t=\chi_t+\xi_t,
\end{equation}
where $\chi_t$ is the short-term factor in prices and $\xi_t$ is the long-term factor. \cite{Bakker2021MatureProgramming} presented a risk-neutral STLT model in discrete time, which is used in this paper for price modelling. The price is represented by, 

\begin{alignat}{1}
    \tilde{p}_t&=e^{\tilde{\chi}_t+\tilde{\xi}_t},\\
    \tilde{\chi}_t&=\tilde{\chi}_{t-1}-\left(1-e^{-\kappa\Delta t}\right)+\sigma_{\chi}\epsilon_{\chi}\sqrt{\frac{1-e^-2\kappa \Delta t}{2\kappa}},\\
    \tilde{\xi}&=\tilde{\xi}_{t-1}+\mu^{*}_{\xi}\Delta t+\sigma_{\xi}\epsilon_{\xi}\sqrt{\Delta t},
\end{alignat}
where $\tilde{p}_t$, $\tilde{\chi}_{t}$ and $\tilde{\xi}_{t}$ are risk-neutral equivalents to the spot price, short-term factor and long-term factor. The volatilities are represented by $\sigma_{\chi}$ and $\sigma_{\xi}$, while $\epsilon_{\chi}$ and $\epsilon_{\xi}$ are correlated standard normal random variables with correlation coefficient $\rho_{\chi \xi}$. The parameter $\kappa$ is the mean reversion coefficient, $\lambda_\chi$ is a risk premium that specifies a reduction in the drift got the short-term process, and $\mu^{*}_{\xi}$ is the risk-neutral drift of the equilibrium level, $\tilde{\xi}_{t}$. The length of the time period $t$ (year) is represented by $\Delta t$.

Based on the historical record, we assume the gas price is closely correlated with the oil price, and so we use a single price process for oil and gas prices.

Then we obtain the mean values of the prices for each stage and construct the mean scenario. Next, we use a Julia package \texttt{ScenTrees.jl} \citep{Kirui2020ScenTrees.jl:Programming} that utilises the methodology proposed in \cite{Pflug2015DynamicTrees} to generate a multi-stage scenario tree. The structure of the scenario tree is not determined a priori but dynamically adapted to meet the Wasserstein distance requirement, which measures of the quality of the approximation. Note that because we consider long-term time horizon with a five year gap between stages, the temporal correlation is weak. However, due to the use of Wasserstein distance in the scenario generation routine, the temporal correlation is still preserved. In addition, there is a perfect spatial correlation because the oil and gas prices are assumed to be the same for all regions.

\subsubsection{Geographical representation of the problem}
The problem potentially consists of many regions and results in a large model. Therefore, we use a k-means cluster approach to group platform locations into representative locations, adapted from \cite{zhang2022_OEH}. We use dedicated locations to represent the regions that require higher resolution. All generators and storage units within each aggregated region that have the same characteristics are aggregated into clusters. In this way, within aggregated regions, the model does not make the investment in individual units but in the total for that technology. For retrofit and abandonment decisions on the offshore platforms, we aggregate the platforms in each offshore region and use the representative platforms in the model.

\subsubsection{Modelling assumptions}
We assume that: (a) the Kirchhoff voltage law is omitted and we use a linear direct current power flow model for the power system part, (b) the initial storage level of storage units in each operational node are assumed to be half of their capacities, (c) the switch from natural gas to hydrogen has little impact on the capacity of a pipeline to transport energy \citep{siemens_2021}, (d) linepack in hydrogen pipelines is omitted, (e) investment in new wells is not considered, (f) we simplify the modelling of pressure and temperature of the production processes on platforms and use typical values from the North Sea region, (g) there is no more oil and gas profit once a platform is retrofitted, and the gas profit, associated with a pipeline is lost once it is retrofitted.

\section{Benders decomposition in MHSP}
\label{sec:Benders decomposition in MHSP}
MHSP has a structure that allows the application of Benders-type decomposition for solving large-scale stochastic programming problems. In the following, we first present how a MHSP problem can be decomposed by Benders decomposition. We then propose a stabilised adaptive Benders decomposition for solving the proposed model efficiently. 
\subsection{Benders decomposition in MHSP}
\label{sec:Benders decomposition in general MHSP}
Here, we give a general MHSP formulation of our combined strategic and operational planning problem and show that it has a block separable structure which allows it to be decomposed by Benders decomposition into a single master strategic investment problem and a family of operational problems.

The MHSP formulation is as follows:

\begin{subequations}
\label{eq:MHSP}
\begin{alignat}{3}
    & \min_{\mathbf{x} \in \mathcal{X}}  &&\sum_{i \in \mathcal{I}} \pi_{i}\left(c_{i}^{\top}x_{i} + \min_{y_i \in \mathcal{Y}}q_{i}^{\top}Qy_{i}\right)\\
    & \text{s.t.} && T^{I}_{A_i}x_{A_i} + W^{I}_i x_i \leq h^{I}_i, \phantom{11111} & i \in \mathcal{I} \setminus \{1\}, \label{eq:MHSP_1}\\
    & \quad && W^{I}_1 x_1 \leq h^{I}_1, & \label{eq:MHSP_2}\\
    & \quad && T^O x^{O}_i + W^O y_{i} \leq h^{O}, & i \in \mathcal{I},\label{eq:MHSP_3}
\end{alignat}
\end{subequations}
where $\mathcal{I}$ is the set of strategic nodes and $A_i$ is the ancestor node of node $i$, $x_i$ is the vector of variables for the strategic decisions and the investment state at investment node $i$, and $\mathbf{x}$ is the collection of all the $x_i, i \in \mathcal{I}$. The $\pi_i$ is the probability of strategic node $i$, and the sum of $\pi_i$ in each strategic stage is equal to 1. $\mathcal{X}$ represents the domains of the elements of $\mathbf{x}$, i.e. their bounds  and whether they are integer or continuous. $y_i$ is the vector of operational variables at node $i$, and $\mathcal{Y}$ represents the bounds on the elements of $y_i$, all of which are continuous variables. Together $\mathcal{X}, \mathcal{Y}$ and Constraints \eqref{eq:MHSP_1}-\eqref{eq:MHSP_3} define the feasible region. $x_i^O$ is the part of $x_i$ that affects the subproblem $i$. The dimensions of $x^O_i, q_i$ and $y_i$ are the same in each subproblem $i$ but the dimensions of $x_i$ may differ. The dimension of the $x_i^O$ and of the $q_i$ are typically much smaller. The vectors, variables, and matrices have compatible dimensions. Some elements of $x_i^O$ may be fixed to exogenous parameters, which is equivalent to having different right hand sides in Constraint \eqref{eq:MHSP_3}.

In MHSP, the complicating variables are the strategic decisions, $\mathbf{x}$, that link all the decision nodes. By fixing the complicating variables $\mathbf{x}$, we can decompose the full size problem using Benders-type decomposition. For a given node $i$, the subproblem is formulated as

\begin{subequations}
\begin{alignat}{3}
    g(x^O_i,q_i):= & \min_{y_i \in \mathcal{Y}}  && q_{i}^{\top}Qy_{i}\\
    & \text{s.t.} && T^Ox^O_i + W^Oy_{i} \leq h^{O}, \qquad i \in \mathcal{I}.
\end{alignat}
\label{eq:benders_subproblem_general}
\end{subequations}
The single subproblem $g(x^{O}_i,q_i)$ can include multiple operational scenarios. If the operational scenarios are independent of each other then an alternative approach would be to treat each as an independent subproblem. This has the potential advantage of maintaining a more accurate model of the operational problems, but the disadvantage is that it increases the number of Benders cuts and slows down the solution of each Benders reduced master problem. In this paper, the single subproblem approach is taken, as illustrated in Figure \ref{fig:MHSP}. The Reduced Master Problem (RMP) at iteration $j$ is 

\begin{subequations}
\begin{alignat}{3}
    & \min_{\mathbf{x \in \mathcal{X}}}  &&\sum_{i \in \mathcal{I}} \pi_{i}(c_{i}^{\top}x_{i} + \beta_i)\\
    & \text{s.t.} && \mathrlap{\text{Constraints }\eqref{eq:MHSP_1} \text{ and } \eqref{eq:MHSP_2},}\notag\\
    &\quad  && \beta_{i} \geq \theta_{ik} + \mathbf{\lambda}^{\top}_{ik}(x^O_i-\hat{x}^{O}_{ik}),
        & \hspace{2cm} i \in \mathcal{I}, k<j, \label{eq:benders_cut_general}
\end{alignat}
\label{eq:benders_RMP_general}
\end{subequations}
where $\beta_i$ is a variable for the approximated cost of the operational subproblem that is embedded in strategic node $i$. Constraint \eqref{eq:benders_cut_general} are the Benders cuts associated with subproblem $i$ built in the iterations up to the current iteration $j$.

\subsection{Stabilised adaptive Benders decomposition}
\label{sec:enhanced benders}

In this section, we propose a stabilised adaptive Benders decomposition algorithm that utilises the structure of MHSP.
The algorithm is taken from \cite{zhang2022_benders} and extended to solve MILP. \cite{zhang2022_benders} utilised the adaptive oracles proposed in \cite{Mazzi2020} and level method stabilisation and achieved a significant reduction in solution time compared to the unstabilised version of adaptive Benders. The adaptive oracles were introduced for problems where the following conditions hold:

\begin{cond}
    $g(x^{O}_i,q_i)$ is convex w.r.t. the vector $x^{O}_i$, and concave w.r.t. the vector $q_i$, and $g(x^{O}_i,q_i)$ is a decreasing function of the elements of $x^{O}_i$ and an increasing function of the elements of $q_i$.  \label{propty: g(x,q)}
\end{cond}

The convexity and concavity are immediate consequences of Equations \eqref{eq:benders_subproblem_general} being a minimisation linear program and the monotonicity properties hold if, for example, $T^O, W^{O}$ and $y_i$ are non-negative. 

Once one or more subproblems have been solved at a collection of points, this information can be used by the adaptive oracles to generate valid bounds at different points for all subproblems. We refer to \cite{Mazzi2020} for the mathematical definition and the proof of the properties of the adaptive oracles.

The adaptive oracles provide bounds for a subproblem at a new solution point without having to solve it exactly, and this reduces the computational cost compared to standard Benders decomposition. The process is shown in Algorithm \ref{alg:Level set method stabilised Benders decomposition with adaptive oracles}. In iteration \(j\), when subproblem \(\hat{i}\) is solved at the point \((\hat{x}^O_{\hat{i}j}, q_{\hat{i}})\), the solver returns the optimal value \(\theta_{\hat{i}j} = g(\hat{x}^O_{\hat{i}j}, q_{\hat{i}})\), and the subgradients \(\lambda_{\hat{i}j}\) and \(\phi_{\hat{i}j}\) with respect to \(x^O_{\hat{i}j}\) and \(q_{\hat{i}}\), respectively. Then the solution vector \((\hat{x}^O_{\hat{i}j}, q_{\hat{i}}, \theta_{\hat{i}j}, \lambda_{\hat{i}j}, \phi_{\hat{i}j})\) is added to the collection \(\mathcal{S}\) of solution vectors. Then, using the information in \(\mathcal{S}\), the adaptive oracles generate valid bounds for all subproblems: the oracles are called for each subproblem \(i\) at the current solution point \((\hat{x}^O_{ij}, q_i)\) and return the values \(\underline{\theta}_{ij}, \overline{\theta}_{ij}\), and \(\lambda_{ij}\) with the properties:
\begin{proper}$\underline{\theta}_{ij} + \lambda^T_{ij}(x^O_i - \hat{x}^O_{ij}) \le g(x^O_i, q_i),~~ \forall x^O_i~~~ \text{and} ~~~ g(\hat{x}^O_{ij}, q_j) \le \overline{\theta}_{ij}.
$ \label{propty: cuts}
\end{proper}

\begin{figure}[!htb]
    \centering
    \includegraphics[scale=0.5]{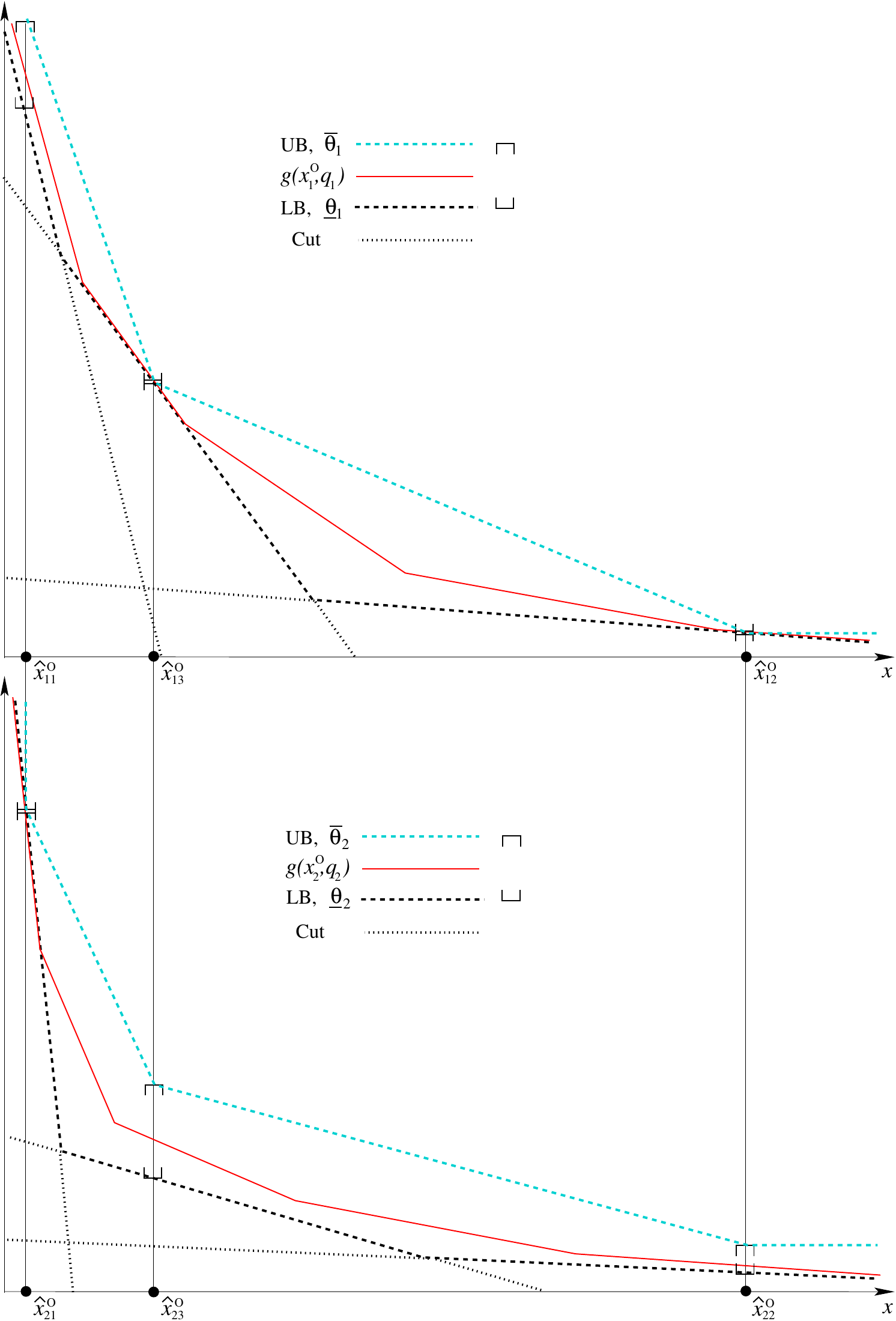}
    \caption{Illustration of adaptive oracles. There are three iterations in this illustration. The \(x\)-axis in the graph shows the value of the first element of \(\hat{x}^O_{ij}\). In the first iteration, subproblem 2 is solved exactly at \((\hat{x}^O_{21}, q_2)\), and an exact cut is generated. An inexact but valid cut and upper bound are then generated by the adaptive oracles for subproblem 1 via cut sharing. In iterations 2 and 3, the exact cut is for subproblem 1, and valid inexact cuts are generated for subproblem 2.}
    \label{fig:algorithm_illustration}
\end{figure}

The RMP in adaptive Benders is the same as in standard Benders, except that the exact cuts in Equation \eqref{eq:benders_cut_general} of standard Benders are replaced by the approximate cuts in Equation \eqref{eq:inexact_cuts}, which use the quantities supplied by the adaptive oracles. An illustration of cut sharing using adaptive oracles is presented in Figure \ref{fig:algorithm_illustration}.

In the standard unstabilised version of Benders the subproblems are evaluated at the current RMP solution. In stabilised versions, they are evaluated instead at a reference point. In the level method stabilisation approach shown in Equations \eqref{eq:LMP}, the reference point, $\mathbf{x}^{Ref}_j$, for iteration $j$ is chosen as the point at the minimum distance from the previous reference point, $\mathbf{x}^{Ref}_{j-1}$, subject to being feasible for the RMP and having a RMP objective value not exceeding the specified target level, $T_j$. The target level $T_j$ is defined in terms of a stabilisation factor $\gamma_j$ by $T_j=L^*_j+\gamma_j(U^*_{j-1}-L^*_j)$, where $U^*_{j-1}$ and $L^*_j$ are the best upper bound and lower bound respectively at iterations $j-1$ and $j$.

\begin{subequations}
    \begin{alignat}{3}
&\mathbf{x}^{Ref}_j:=\argmin_{\mathbf{x}\in \mathcal{X},\beta} &&\norm{\mathbf{x}-\mathbf{x}^{Ref}_{j-1}}^2_2\\
        & \text{s.t.} && \mathrlap{\text{Constraints }\eqref{eq:MHSP_1} \text{ and } \eqref{eq:MHSP_2},}\notag\\
        &\quad  && \beta_{i} \geq \underline{\theta}_{ik} + \underline{\mathbf{\lambda}}^{\top}_{ik}(x^O_i-\hat{x}^{O}_{ik}),
        & \hspace{1cm} i \in \mathcal{I}, k<j, \label{eq:inexact_cuts}\\
        &\quad && \sum_{i \in \mathcal{I}} \pi_i(c^\top_i x_i + \beta_i) \le T_j, \label{eq:level set}
    \end{alignat}
    \label{eq:LMP}
\end{subequations}

The level method stabilisation problem is a quadratic program if $L_2$ norm is used \citep{Zverovich2012} and becomes a mixed-integer quadratic program when integer variables are present. In this paper, we instead use a centred point approach to avoid solving the mixed-integer quadratic program but still obtain a stabilised solution. 

\begin{figure}[!htb]
    \centering
    \includegraphics[scale=0.4]{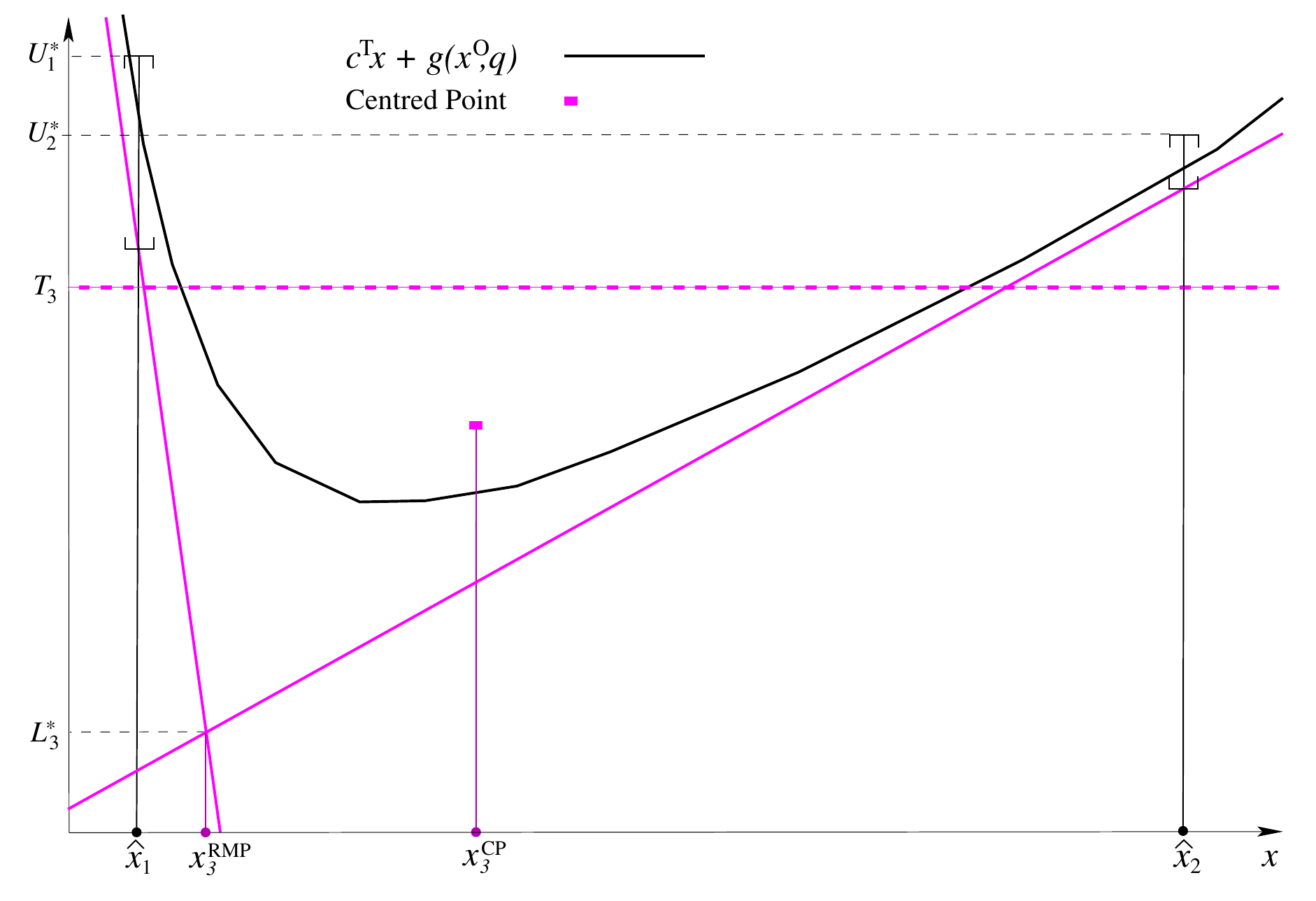}
    \caption{Illustration of centred point stabilisation. We assume $x$ has only one element and the inexact oracles provide inexact but valid cuts and upper bounds. The centred point is relative to the two inexact cuts and the target $T_3$.}
    \label{fig:centred_point}
\end{figure}

The Centred Point problem (CP) is derived from the level method stabilisation problem by relaxing the integer variables and dropping the objective function.  The CP is, therefore, a linear programming feasibility problem.  An illustration of the centred point stabilisation approach is presented in Figure \ref{fig:centred_point}. Finding a centred point in the feasible region can be done efficiently using solvers like Gurobi \citep{gurobi}. For example, using the barrier algorithm of Gurobi and turning off presolve and crossover, the CP finds a well centred point within the feasible region approximating the analytic centre, (i.e. the point that maximises the product of distances from all faces of the linear programming polytope \citep{Gondzio1996}).  Solving a CP is computationally much cheaper than solving the level method stabilisation problem, so the time spent on stabilisation is reduced significantly.

The method is presented in Algorithm \ref{alg:Level set method stabilised Benders decomposition with adaptive oracles}.  The stabilisation factor $\gamma$ is adjusted using an approach that is analogous to the trust region method. In \cite{zhang2022_benders}, it was shown that by adjusting the stabilisation factor, a Benders type algorithm became more robust. In the early tests in this paper, we noticed that adjusted stabilisation is much more efficient than using a fixed stabilisation factor.

The level set adjustment steps are presented in Subroutine \ref{alg:level set management}. If both the actual improvement and predicted improvement are positive (Line \ref{alg:positive_ratio}), we check the ratio against two predefined parameters $\underline{P}$ and $\overline{P}$. If the improvement ratio $r$ is lower than $\underline{P}$, meaning that the actual improvement is less than the least improvement we want to achieve in proportion to the predicted improvement, we tighten the stabilisation in Line \ref{alg:tighten}. This is because we want to avoid moving to an area that provides insufficient improvement. The parameter $\omega$ affects how significant the change in stabilisation is. If the improvement ratio is higher than $\overline{P}$, the actual improvement aligns with the predicted improvement well. This suggests we are moving in the right direction, and we loosen the stabilisation to hopefully achieve more improvement in the next iteration.

\begin{algorithm}[!htb]
    \caption{Stabilised Benders decomposition with adaptive oracles}\label{alg:Level set method stabilised Benders decomposition with adaptive oracles}
    \begin{spacing}{0.9}
    \begin{algorithmic}[1]
        \State choose $\epsilon$ (convergence tolerance), $\gamma_0 \in (0,1)$ (stabilisation factor), $\underline{\beta}$ (initial lower bound $\beta_{i}$), $U^{*}_0:=+\infty$ (initial upper bound), $\epsilon^{0}$, $\omega \in (0,1)$, $\underline{P} \in (0,1)$, $\overline{P} \in (\underline{P},1)$, $\xi:=0, \mathbf{x}^{CP}_i:=0$;
        \State solve subproblem exactly at the special point $(\underline{x}^{O},\underline{q})$ and obtain $\theta$, $\lambda$ and $\phi$; $\mathcal{S}:=\{(\underline{x}^O,\underline{q},\theta,\lambda,\phi)\}$;
        \Repeat \label{alg: outer start}
            \State $j:=j+1$;
            \State solve RMP and obtain $\beta_{ij}$ and $\mathbf{x}^{RMP}_{j}$; $L^{*}_{j}:=\sum_{i \in \mathcal{I}}\pi_{i}(c^{\top}x^{RMP}_{ij} +\beta_{ij})$;\label{alg:solve RMP}
        \StatePar{$\mathbf{\hat{x}}_{j}:=\mathbf{x}^{RMP}_{j}$;}
            \If{$j \geq 2$ \textbf{and} $\xi=1$}{}
            \StatePar{set CP target $T_j:=L^*_j+\gamma_j(U^*_{j-1}-L^*_j)$, solve CP and obtain $\mathbf{x}^{CP}_j$;}
            \If{$\norm{\mathbf{x}^{CP}_{j}-\mathbf{x}^{CP}_{j-1}} \leq \epsilon^{0}$}{ $\mathbf{\hat{x}}_{j}:=\mathbf{x}^{CP}_{j}$;}
            \EndIf
            \Else
            \StatePar{$\mathbf{x}^{CP}_j:=\mathbf{x}^{CP}_{j-1}$;}
            \EndIf
            \For{$i \in \mathcal{I}$}
                \StatePar{call adaptive oracles at $(\hat{x}^O_{ij},q_{i})$ and  obtain $\underline{\theta}_{ij}$, $\overline{\theta}_{ij}$, and $\underline{\lambda}_{ij}$;}
            \EndFor
            \StatePar{set $\xi:=0$;}
            \Repeat
               \If{$\max \pi_i(\overline{\theta}_{ij}-\underline{\theta}_{ij})=0$}{ Exit loop}\EndIf \label{alg: check gap}
        \State $\hat{i}:= \argmax_{i \in \mathcal{I}} \pi_i (\overline{\theta}_{ij} - \underline{\theta}_{ij})$; 
        \State solve subproblem exactly at $(\hat{x}^O_{\hat{i}j}, q_{\hat{i}})$ and obtain $\theta_{\hat{i}j}$, $\lambda_{\hat{i}j}$, $\phi_{\hat{i}j}$;
        \If{$\theta_{\hat{i}j}>\underline{\theta}_{\hat{i}j}$}{}
        \StatePar{set $\xi:=1$;}
        \State $\mathcal{S}:=\mathcal{S}\cup\{(\hat{x}^O_{\hat{i}j}, q_{\hat{i}},\theta_{\hat{i}j},\lambda_{\hat{i}j},\phi_{\hat{i}j})\}$;
        \EndIf
                \For{$i \in \mathcal{I}$}
                    \StatePar{call adaptive oracles at $(\hat{x}^{O}_{ij},q_{i})$ and obtain $\underline{\theta}_{ij}$, $\overline{\theta}_{ij}$, and $\underline{\lambda}_{ij}$;}
                \EndFor
                \StatePar{$\hat{L}_{j}:= \sum_{i \in \mathcal{I}}\pi_{i}(c^{\top}\hat{x}_{ij}+\underline{\theta}_{ij})$;}
                \StatePar{$\hat{U}_{j}:= \sum_{i \in \mathcal{I}}\pi_{i}(c^{\top}\hat{x}_{ij}+\overline{\theta}_{ij})$;}
            \Until{$(\hat{U}_{j}-\hat{L}_{j} \leq U^{*}_{j-1}-L^{*}_{j-1} \text{ and } \xi=1)$ \normalfont{\textbf{or}} $\hat{L}_{j}\geq U^{*}_{j-1}$} \label{alg: exit inner loop}
            \For{$i \in \mathcal{I}$}
            \StatePar{call upper bound adaptive oracle at $(x_{ij}^{RMP,O}, q_i)$ and obtain $\overline{\theta}_{ij}$;}
            \EndFor
            \StatePar{$U_{j}:=\sum_{i \in \mathcal{I}}\pi_{i}(c^{\top}x^{RMP}_{ij} +\overline{\theta}_{ij})$, $U^{*}_{j}:=\min(U^{*}_{j-1},U_{j})$; \label{alg:update upper bound}}
            \StatePar{Subroutine \ref{alg:level set management}: Dynamically adjust the level set;}
        \Until{$U^{*}_{j}-L^{*}_{j} \leq \epsilon.$} \label{alg:outer end}
    \end{algorithmic}
     \end{spacing}
\end{algorithm}

\begin{subroutine}[!htb]
    \caption{Dynamic level set management}\label{alg:level set management}
    \begin{spacing}{0.9}
    \begin{algorithmic}[1]
        \StatePar{$I^{P}_{j}:=\hat{L}_{j-1}-{T_{j}}$, $I^{A}_{j}:=\hat{L}_{j-1}-\hat{L}_{j}$, $r:=\frac{I^{A}_{j}}{I^{P}_{j}}$;}
        \If{$I^{A}_{j}>0$ \textbf{and} $I^{P}_{j}>0$}{} \label{alg:positive_ratio}
        \If{$r \leq \underline{P}$}{}
        \State{$\gamma_{j+1}:=1-\omega(1-\gamma_j)$;} \label{alg:tighten}
        \ElsIf{$r>\overline{P}$}{}
        \State{$\gamma_{j+1}:=\omega\gamma_j$;} \label{alg:loosen}
        \EndIf
        \EndIf
    \end{algorithmic}
    \end{spacing}
\end{subroutine}

\subsection{Convergence of Algorithm \ref{alg:Level set method stabilised Benders decomposition with adaptive oracles}}
In this section, we prove the convergence of Algorithm \ref{alg:Level set method stabilised Benders decomposition with adaptive oracles}.

\begin{lemma}
For problems where Condition \ref{propty: g(x,q)} holds, Algorithm \ref{alg:Level set method stabilised Benders decomposition with adaptive oracles} achieves an $\epsilon$-optimal solution when it terminates.
\end{lemma}

\begin{proof}
It is proved in \cite{Mazzi2020} that if Condition \ref{propty: g(x,q)} holds, then the bounds in Property \ref{propty: cuts} provided by the adaptive oracles are valid. It then follows that the RMP, Equations \eqref{eq:benders_RMP_general}, is a valid relaxation of the MHSP, Equations \eqref{eq:MHSP}, and consequently that \(L^*_j\) in Line \ref{alg:solve RMP} is a valid lower bound on the optimal MHSP objective.

Also, \(U_j\) in Line \ref{alg:update upper bound} is obtained from an integer feasible value of \(x^{RMP}_{ij}\) and an upper bound on the objective value for the subproblems evaluated at the RMP point. Therefore, \(U_j\) is a valid upper bound for the problem, and so is \(U^*_j\).

Since \(L^*_j\) and \(U^*_j\) are valid upper and lower bounds on the optimal objective, the algorithm will terminate correctly if the convergence tolerance is met.
\end{proof}

\begin{theorem}
If all the subproblems are linear programming problems, then Algorithm \ref{alg:Level set method stabilised Benders decomposition with adaptive oracles} converges in a finite number of iterations to an $\epsilon$-optimal solution.
\end{theorem}

\begin{proof}
When a subproblem is solved, the cut generated is a face of its value surface. If the new bound at the point is tighter than before, i.e., \(\theta_{ij} > \underline{\theta}_{ij}\), then the cut is new.

If no new face is generated in the last outer iteration (Lines \ref{alg: outer start}-\ref{alg:outer end}), then the next iteration is an RMP iteration. Hence, the number of CP iterations must be less than or equal to the total number of faces, which is finite because the subproblems are linear programs.

Also, there can only be a finite number of RMP iterations at which a new face is generated. Eventually, if termination has not already occurred, there must be an RMP iteration where no new face is found. When that happens, the inner loop exits either at Line \ref{alg: check gap}, or with the \(\hat{L}_j \ge U^*_{j-1}\) condition at Line \ref{alg: exit inner loop}. In the first case, there is no bound gap at the optimal RMP point. In the second case, since this is an RMP iteration, \(\hat{L}_j\) has been evaluated at the RMP point and is therefore a lower bound on the RMP, and since this is greater than or equal to the upper bound \(U^*_{j-1}\), equality must hold. Hence, in either case, the exact optimum has been reached, so the termination condition has been satisfied within a finite number of outer iterations.

Finally, note that when the adaptive oracles are called at a point where subproblem \(\hat{i}\) has been previously solved, they return a zero bound gap, i.e., \(\overline{\theta}_{\hat{ij}} - \underline{\theta}_{\hat{ij}} = 0\). Also, as more exact solutions are added to $\mathcal{S}$, the oracle bounds cannot become slacker, so the bound gap at a point where the subproblem has previously been solved remains zero. Hence, the inner loop cannot select a subproblem that has been solved in an earlier inner iteration, so the loop cannot be executed more than $|\mathcal{I}|$ times. Hence, since the number of outer iterations is also finite, the total number of function evaluations is finite.
\end{proof}

\section{The REORIENT model}
\label{sec:model}
This section presents the energy system integrated planning and operational optimisation model. The full model is decomposed by having an investment planning master problem and an operational subproblem. We use the conventions that calligraphic capitalised Roman letters denote sets, upper case Roman and lower case Greek letters denote parameters, and lower case Roman letters denote variables. The indices are subscripts, and name extensions are superscripts. The same lead symbol represents the same type of thing. The names of variables, parameters, sets and indices are single symbols. We give a brief definition of some of the main sets and variables, and their corresponding domains as we explain the equations. For a complete overview of all sets and indices, parameters and variables used in the REORIENT model, we refer to \ref{sec:nonmenclature}.

\subsection{Investment planning model (RMP in Benders decomposition)}
The investment master problem Equations \eqref{eq:investment_cost}-\eqref{eq:mp_domain} follows the general formulation given by Equations \eqref{eq:benders_RMP_general}. The total discounted cost for investment planning, Equation \eqref{eq:investment_cost}, consists of actual investment costs $c^{INV}$ as well as the expected operational cost of the system over the time horizon $\kappa\sum_{i \in \mathcal{I}^{Ope}}\pi_{i}c^{OPE}_i$ which is total approximated subproblem costs in Benders decomposition. Here, $\kappa$ is a scaling factor that depends on the time step between two successive investment nodes. The scaling factor scales the operational costs between two successive investment nodes. By doing this, we can evaluate the operational subproblem on the represented operational hours and scale the cost up to obtain the total operational costs. Equation \eqref{eq:cinv} calculates the investment cost, which comprises capacity-dependent and capacity-independent costs, retrofitting costs, abandonment costs, fixed operating and maintenance costs, and profit of existing technology (e.g., oil and gas platform).

\begin{equation}
    \min c^{INV}+\kappa\sum_{i \in \mathcal{I}^{Ope}}\pi_{i}c^{OPE}_i ,\label{eq:investment_cost}
\end{equation}
where

\begin{equation}
\begin{split}
    c^{INV}=\sum_{i\in \mathcal{I}^{Inv}}\pi^{Inv}_{i}\left(\sum_{p \in \mathcal{P}}\left(C^{InvV}_{pi}x^{Inv}_{pi}+C^{InvF}_{pi}y^{Inv}_{pi}\right)+\sum_{p \in \mathcal{P}^{RT}}\left(C^{ReTV}_{pi}x^{ReT}_{pi}+C^{ReTF}_{pi}y^{ReT}_{pi}\right)\right)+\\
    \phantom{ass}\kappa \sum_{i \in \mathcal{I}^{Ope}}\pi^{Ope}_{i}\left(\sum_{p \in \mathcal{P}}C^{Fix}_{pi}x^{Acc}_{pi} + \sum_{p \in \mathcal{P}^{RT}}C^{ReTFixO}_{pi}x^{AccReT}_{pi}+\sum_{p \in \mathcal{P}^{R}}C^{ReFFixO}_{pi}x^{AccReF}_{pi} \right). \label{eq:cinv}
\end{split}
\end{equation}

Constraint \eqref{eq:cap_tech} states that the accumulated capacity of a technology $x^{Acc}_{pi}$ in an operational node equals the sum of the historical capacity $X^{Hist}_{p}$ and newly invested capacities $x^{Inv}_{pi}$ in its ancestor investment nodes $\mathcal{I}^{Inv}_{i}$ that are not retired.

\begin{equation}
    x^{Acc}_{pi}=X^{Hist}_{pi}+\sum_{j \in\mathcal{I}^{Inv}_{i}|\kappa(S^{Ope}_{i}-S^{Inv}_{j})\leq H_{p}}x_{pj}^{Inv}, \qquad p \in \mathcal{P},  i \in \mathcal{I}^{Ope}.\label{eq:cap_tech}
\end{equation}

Constraint \eqref{eq:max_built} ensures the maximum $X^{MaxInv}_{pi}$ capacity that is built in an investment node. The binary variable $y^{Inv}_{pi}$ equals 1 if a technology $p\in \mathcal{P}$ in investment node $i \in \mathcal{I}^{Inv}$ is invested. Parameter $X^{MaxAcc}_{p}$ gives the maximum capacity that can be installed for different technologies.

\begin{equation}
    x^{Inv}_{pi} \leq X^{MaxInv}_{pi}y^{Inv}_{pi}, \qquad p \in \mathcal{P}, i \in \mathcal{I}^{Inv}. \label{eq:max_built}
\end{equation}

Constraints \eqref{eq:max_acc}-\eqref{eq:retrofit_cap} establish that the invested capacity and accumulated capacity of newly invested technologies and retrofitted technologies should be within the capacity limits. 

\begin{equation}
    x^{Acc}_{pi} \leq X^{MaxAcc}_{p},  \qquad  p \in \mathcal{P},  i \in \mathcal{I}^{Ope}, \label{eq:max_acc}
\end{equation}
\begin{equation}
    x^{AccReT}_{pi} \leq X^{MaxAccReT}_{p}, \qquad  p \in \mathcal{P}^{RT},  i \in \mathcal{I}^{Ope}, \label{eq:max_acc_retrofit}
\end{equation}
\begin{equation}
    x^{ReT}_{pj} \leq X^{MaxReT}_{p}y^{ReT}_{pi}, \qquad p \in \mathcal{P}^{RT}, i \in \mathcal{I}^{Inv},j \in \mathcal{I}^{Inv}_{i}. \label{eq:retrofit_cap}
\end{equation}

Constraint \eqref{eq:retrofit_from} dictates that the existing capacity is zero if a technology is retrofitted to a new technology. 

\begin{equation}
    x^{AccReF}_{pj} = X^{HistReF}_{pj}(1-y^{ReF}_{pi}), \qquad p \in \mathcal{P}^{R}, i \in \mathcal{I}^{Inv}, j \in \mathcal{I}^{Ope}_{i}. \label{eq:retrofit_from}
\end{equation}

Constraint \eqref{eq:retrofit_to_limits} states that only one technology can be retrofitted to. 

\begin{equation}
    \sum_{p \in \mathcal{P}^{R}_{p}}y^{ReT}_{pi} \leq y^{ReF}_{pi}, \qquad p \in \mathcal{P}^{R}, i \in \mathcal{I}^{Inv}. \label{eq:retrofit_to_limits}
\end{equation}

Constraint \eqref{eq:one_time_retrofit} ensures that retrofit can only happen once for a technology during the planning horizon. 

\begin{equation}
    \sum_{i \in \mathcal{I}^{Inv}}y^{ReF}_{pi} \leq 1, \qquad p \in \mathcal{P}^{R}. \label{eq:one_time_retrofit}
\end{equation}

Constraint \eqref{eq:retrofit_acc} states that the accumulated capacity of a technology $x^{AccReT}_{pi}$ in an operational node equals the newly invested capacities $x^{ReT}_{pi}$ in its ancestor investment nodes $\mathcal{I}^{Inv}_{i}$ that are not retired.  Parameter $X^{MaxAccReT}_{p}$ is the maximum accumulated capacity of a technology that is retrofitted from another. 

\begin{equation}
    x^{AccReT}_{pi} = \sum_{j \in\mathcal{I}^{Inv}_{i}|\kappa(S^{Ope}_{i}-S^{Inv}_{j})\leq H_{p}}x^{ReT}_{pj}, \qquad  p \in \mathcal{P}^{RT}, i \in \mathcal{I}^{Ope}. \label{eq:retrofit_acc}
\end{equation}

The Benders cuts built up to iteration $k-1$ are given by Equation \eqref{eq:benders_cut}.

\begin{equation}
    c^{OPE}_i \geq \theta + \lambda^{\top}(x_{i}-x), \qquad (x,\theta,\lambda) \in \mathcal{F}_{i(k-1)}, i \in \mathcal{I}. \label{eq:benders_cut}
\end{equation}

The domains of variables are given as follows

\begin{equation}
   x^{Inv}_{pi}, x^{Acc}_{pi}, x^{ReT}_{pi}, x^{AccReT}_{pi}, x^{AccReF}_{pi}, c^{INV} \in \mathbb{R}^{+}_{0}, \qquad y^{Inv}_{pi}, y^{ReF}_{pi}, y^{ReT}_{pi} \in \{0,1\}. \label{eq:mp_domain}
\end{equation}

The vector $x_i=\left( \{x^{Acc}_{pi}, p \in \mathcal{P}\}, \{x^{AccReT}_{pi}, p \in \mathcal{P}^{RT}\}, \{x^{AccReF}_{pi}, p \in \mathcal{P}^{R}\}, \mu^{DP}_i, \mu^{DH}_i, \mu^{DHy}_i, \mu^{E}_i\right)$ collects all right-hand side coefficients that will be fixed in operational subproblem, Equations \eqref{eq:SP_objective}-\eqref{eq:sp_domain}. The vector $c_i=\left(C^{CO_2}_{i}\right)$ collects all the cost coefficients. The vectors $x_i$ and $c_i$ will be fixed as parameters in the operational problem. The long-term uncertain parameters,  including load scaling $\mu^{DP}_i$, $\mu^{DH}_i$, and $\mu^{DHy}_i$ $\mu^E_i$, are fixed in the operational problem described below because they affect the system operation.

\subsection{Operational problem (subproblem in Benders decomposition)}
We now present the operational problem and note that we omit index $i$ in the operational model for ease of notation because all variables and parameters are defined for each operational node.

The right hand side parameters $P^{Acc}_p$, $V^{Acc}_c$, $P^{AccG}_g$, $P^{AccHRor}_g$, $P^{AccSE}$, $Q^{AccSE}_s$, $P^{AccL}_l$, $V^{AccLHy}_l$, $\mu^{DP}$, $\mu^{DH}$, $\mu^E$, $\mu^{DHy}$ are fixed by the solution $x_i$ from solving the master problem Equations \eqref{eq:investment_cost}-\eqref{eq:mp_domain}. The CO$_2$ cost of generators that is included in parameter $C^G_g$ is fixed by $c_i$ from the master problem.

The operational cost $c^{OPE}(x_i,c_i)$ at one operational node $i$ is computed by solving subproblem Equations \eqref{eq:SP_objective}-\eqref{eq:sp_domain} given the decisions $x_i$ and $c_i$ made in Equations \eqref{eq:investment_cost}-\eqref{eq:benders_cut}. The operational subproblem Equations \eqref{eq:SP_objective}-\eqref{eq:sp_domain} correspond to the general formulation Equations \eqref{eq:benders_subproblem_general}. The objective function, the operational cost, includes total operating costs of generators $C^{G}_{g}p^{G}_{gt}$, energy load shedding costs for heat, power, and hydrogen $C^{Shed}p_{zt}^{Shed}$ and $C^{Shed}v_{zt}^{Shed}$ and fuel cost of steam reforming plants $C^{R}v^{R}_{rt}$. $C^{G}_{g}$ includes the variable operational cost, fuel cost and the CO$_2$ tax, $C^{CO_2}$, charged on the emissions of generator $g$. The inclusion of load shedding variables $p^{Shed}_{zt}$ and $v^{ShedHy}$ ensures the operational problem is always feasible. The load shedding costs $C^{Shed}$ are large enough so that the optimal solution has little or nor load shed.

\begin{equation}
    \min \sum_{t \in \mathcal{T}}\pi_{t}H_t\left(\sum_{g \in \mathcal{G}}C_g^{G}p_{gt}^{G}+\sum_{r \in \mathcal{R}}C^{R}v_{rt}^{R}+\sum_{z \in \mathcal{Z}}\left(\sum_{l \in \{H,P\}}C^{Shed,l}p_{zt}^{Shed,l}+C^{ShedHy}v_{zt}^{ShedHy}\right)\right). \label{eq:SP_objective}
\end{equation}

Constraints \eqref{eq:all_cap} ensure that the technologies operate within their capacity limits.

\begin{subequations}
\label{eq:all_cap}
    \begin{alignat}{3}
    &\quad&&p_{pt}\leq P_{p}^{Acc}, & p \in \mathcal{P}^{*}, t \in \mathcal{T}, \label{eq:multiple_cap}\\
    &\quad&& v_{vt} \leq V_{v}^{Acc},&   v \in \mathcal{V}^{*}, t \in \mathcal{T}, \label{eq:multiple_gas_cap}\\
    &\quad&& p_{gt}^{G}+p_{gt}^{ResG}\leq P_{g}^{AccG}, & g \in \mathcal{G}, t \in \mathcal{T}, \label{eq:gen_cap}\\
    &\quad&& p_{st}^{SE-} + p_{st}^{ResSE} \leq P_{s}^{AccSE}, & s \in \mathcal{S}^{E}, t \in \mathcal{T}, \label{eq:estore_discharge_cap}\\
    &\quad&& q_{st}^{SE} \leq Q_{s}^{AccSE}, & s \in \mathcal{S}^{E}, t \in \mathcal{T}, \label{eq:estore_energy_cap}\\
    &\quad&&-P_{l}^{AccL} \leq p_{lt}^{L} \leq P_{l}^{AccL}, & l \in \mathcal{L}, t \in \mathcal{T}, \label{line_cap}\\
    &\quad&&-V_{l}^{AccLHy} \leq v_{lt}^{LHy} \leq V_{l}^{AccLHy}, & \qquad l \in \mathcal{L}^{Hy}, t \in \mathcal{T}. \label{eq:pipeline_cap}
    \end{alignat}
\end{subequations}

Constraint \eqref{eq:ramp} captures how fast generators can ramp up or ramp down their power output, respectively. 

\begin{equation}
    -\alpha^{G}_{g}P_{g}^{AccG} \leq p_{gt}^{G}+p_{gt}^{ResG}-p_{g(t-1)}^{G}-p_{g(t-1)}^{ResG} \leq \alpha^{G}_{g}P_{g}^{AccG}, \qquad  g \in \mathcal{G}, n \in \mathcal{N}, t \in \mathcal{T}_{n}. \label{eq:ramp}
\end{equation}

Constraint \eqref{eq:reserve} dictates that the spinning reserve of generator $p^{ResG}_{gt}$, plus the reserve of the electricity storage $p^{ResES}_{st}$ must exceed the minimum reserve requirement, where $\sigma^{Res}$ is a percentage of the power load.

\begin{equation}
    \sum_{g \in \mathcal{G}_{z}} p^{ResG}_{gt}+ \sum_{s \in \mathcal{S}^{E}_{z}}p^{ResSE}_{st} \geq \sigma_{z}^{Res}P^{DP}_{zt}, \qquad z \in \mathcal{Z}, t \in \mathcal{T}. \label{eq:reserve}
\end{equation}

Constraints \eqref{eq:hydro_seasonal} and \eqref{eq:hydro_ror} ensure that run-of-the-river hydropower and regulated hydropower production are within their limits and according to the generation profiles, separately.

\begin{equation}
    \sum_{t \in \mathcal{T}_{n}} p_{gt}^{H} \leq \sum_{t \in \mathcal{T}_{n}} P_{gt}^{HSea}, \qquad g \in \mathcal{G}^{HSea}, n \in \mathcal{N}, \label{eq:hydro_seasonal}
\end{equation}
\begin{equation}
    p_{gt}^{H} \leq P_{gt}^{HRor}P^{AccHRor}_g, \qquad g \in \mathcal{G}^{HRor}, t \in \mathcal{T}. \label{eq:hydro_ror}\\ 
\end{equation}

Constraint \eqref{eq:kcl} ensures that, in one operational period $t$, the sum of total power generation of generators $p^{G}_{gt}$, power discharged from all the electricity storage $p^{SE-}_{st}$, renewable generation $R^{R}_{zt}p^{AccR}_{rt}$, hydro power generation $p^{H}_{gt}$, fuel cell generation $p^{F}_{ft}$, power transmitted to this region, and load shed $p^{ShedP}_{zt}$ equals the sum of  power demand $\mu^{DP}P^{DP}_{zt}$, power consumption of electric boilers $p^{BE}_{bt}$, power consumption of all electrolysers $p^{E}_{et}$, power transmitted to other regions, and power generation shed $p^{GShedP}_{zt}$. The parameter $R^{GR}_{rt}$ is the capacity factor of the renewable unit that is a fraction of the nameplate capacity $P^{AccR}$. The subset of a technology in the region $z$ is represented by $R_{z}:=\{r \in \mathcal{R}: r \text{ is available in region } z\}$, where $\mathcal{R}$ can be replaced by other sets of technologies. The power load shed $p^{ShedP}$ allows power demand unmet at a high cost to ensure the feasibility of the operational subproblem. The same idea applies to hydrogen mass balance and heat energy balance.

\begin{equation}
\begin{split}
\sum_{g \in \mathcal{G}_{z}}p_{gt}^{G}+\sum_{l \in \mathcal{L}^{In}_{z}}\eta^{L}p_{lt}^{L}+\sum_{s \in \mathcal{S}^{E}_{z}}p_{st}^{SE-}+\sum_{r \in \mathcal{G}^{R}_{z}}R^{GR}_{rt}P_{r}^{AccGR}+\sum_{g \in \mathcal{G}^{H}_{z}}p_{gt}^{H}+
\sum_{f \in \mathcal{F}_{z}}p^{F}_{ft}+p_{zt}^{ShedP}= \quad\\
 \mu^{DP} P^{DP}_{zt}+\sum_{b \in \mathcal{B}^{E}_{z}}p^{BE}_{bt}+\sum_{e \in \mathcal{E}_{z}}p^{E}_{et}+
\sum_{l \in \mathcal{L}^{Out}_{z}}\eta^{L}_{l}p_{lt}^{L}+\sum_{s \in \mathcal{S}^{E}_{z}}p_{st}^{SE+}+p_{zt}^{GShedP}, \qquad z \in \mathcal{Z}, t \in \mathcal{T}. \label{eq:kcl}
\end{split}
\end{equation}

The hydrogen mass balance Constraint \eqref{hydrogen_balance} dictates that hydrogen produced by electrolyser $H_t\rho^{E}p^{E}_{et}$ and steam reforming plant $v^{R}_{rt}$, hydrogen transmitted to this region, withdraw from a hydrogen storage $v^{SHy-}_{st}$ and hydrogen production shed $v^{GShedHy}_{zt}$ equals the hydrogen demand $V^{DHy}_{zt}$, fuel supply to fuel cell $H_t\rho^{F}p^{F}_{ft}$, hydrogen injected into the storage $v^{SHy+}$, hydrogen transmitted from this region plus the hydrogen load shed $v^{ShedHy}_{zt}$.

\begin{equation}
    \begin{split}
    \sum_{s \in \mathcal{S}^{Hy}_{z}}v^{SHy+}_{st}+\sum_{l \in \mathcal{L}^{HyOut}_{z}}v^{LHy}_{lt}+\sum_{f \in \mathcal{F}_{z}}H_{t}\rho^{F}p^{F}_{ft}+v^{ShedHy}_{zt}+\mu^{DHy}V^{DHy}_{zt}= \phantom{abcdedmsssss}\\
    \sum_{l \in \mathcal{L}^{HyIn}_{z}}v^{LHy}_{lt}+\sum_{e \in \mathcal{E}_{z}}H_t\rho^{E}p^{E}_{et}+\sum_{r \in \mathcal{R}_{z}}v^{R}_{rt}+\sum_{s \in \mathcal{S}^{Hy}_{z}}v^{SHy-}+v^{GShedHy}_{zt}, \qquad z \in \mathcal{Z}, t \in \mathcal{T}.\label{hydrogen_balance}
    \end{split}
\end{equation}

The heat energy balance Constraint \eqref{heat_balance} states that the heat recovery of gas turbines $\eta^{HrG}_{g}p^{G}_{gt}$, plus electric boiler heat generation $\eta^{BE}_{b}p^{BE}_{bt}$, plus heat load shed $p^{ShedH}_{zt}$ equals the heat demand $\mu^{DH}P^{DH}_{zt}$ plus the heat generation shed $p^{GShedH}_{zt}$.

\begin{equation}
    \begin{split}
    \sum_{g \in \mathcal{G}}\eta^{HrG}_gp^{G}_{gt} +\sum_{b\in \mathcal{B}^{E}_{z}}\eta^{BE}_{b}p^{BE}_{bt}+p^{ShedH}_{zt}=\mu^{DH}P^{DH}_{zt}+p^{GShedH}_{zt}, \qquad z \in \mathcal{Z}^{P}, t \in \mathcal{T}. \label{heat_balance}
    \end{split}
\end{equation}

Constraint \eqref{eq:storage_balance} states that the state of charge $q^{SE}_{st}$ in period $t+1$ depends on the previous state of charge $q^{SE}_{st}$, the charged power $\mu^{SE}_sp^{SE+}_{st}$ and discharged power $p^{SE-}_{st}$. The parameter $\eta^{SE}_{s}$ represent the charging efficiency. 

\begin{equation}
    q_{s(t+1)}^{SE}=q_{st}^{SE}+H_{t}(\eta_{s}^{SE}p_{st}^{SE+}-p_{st}^{SE-}), \qquad  s \in \mathcal{S}^{E}, n \in \mathcal{N}, t \in \mathcal{T}_{n}. \label{eq:storage_balance}
\end{equation}

The hydrogen storage balance Constraint \eqref{eq:hydrogen_storage_balance} shows that the hydrogen storage level $v^{SHy}_{st}$ at period $t+1$ equals to storage level at the previous period, plus the hydrogen injected $v^{SHy+}_{st}$, minus the hydrogen withdrawn $v^{SHy-}_{st}$.

\begin{equation}
    v^{SHy}_{s(t+1)}=v^{SHy}_{st}+v^{SHy+}_{st}-v^{SHy-}_{st}, \qquad  s \in \mathcal{S}^{Hy}, n \in \mathcal{N}, t \in \mathcal{T}_{n}. \label{eq:hydrogen_storage_balance}\\
\end{equation}

Constraint \eqref{eq:co2budget} restricts the total emission. The parameter $\mu^{E}$ is the CO$_2$ budget.

\begin{equation}
    \sum_{t \in \mathcal{T}}\pi_{t}\left(\sum_{g \in \mathcal{G}} E^{G}_{g}p^{G}_{gt}+\sum_{r \in \mathcal{R}}E^{R}v^{R}_{rt}\right)\leq \mu^{E}. \label{eq:co2budget}
\end{equation}

The domains of variables are given as follows

\begin{equation}
\begin{split}
    p^{L}_{lt}, v^{LHy}_{lt} \in \mathbb{R},\quad p^{G}_{gt}, p^{ShedP}, p^{ShedH}, v^{ShedHy}_{zt}, p^{Acc}_{p}, p^{BE}_{bt}, p^{ResG}, p^{ResSE}\in \mathbb{R}^{+}_{0}, \\
    p^{AccG}_{g}, p^{GShedP}_{zt}, p^{GShedH}_{zt}, v^{GShedHy}_{zt}, v^{SHy+}_{st}, v^{SHy-}_{st}, v^{SHy}_{st}, p^{E}_{et}, p^{H}_{gt}, p^{AccGR}_{r} \in \mathbb{R}^{+}_{0},\\
    p^{SE+}_{st}, v^{AccLHy}_{l}, p^{SE-}_{st}, q^{AccSE}_{s}, q^{SE}_{st}, p^{AccL}_{l}, p^{AccR}_{r}, p^{F}_{ft}, v^{R}_{rt} \in \mathbb{R}^{+}_{0}.\label{eq:sp_domain}
\end{split}
\end{equation}

\section{Results}
\label{sec:results}
In this section, we first present the case study. Then we report the computational performance of the stabilised adaptive Benders decomposition, followed by the sensitivity analysis of the retrofitting cost of natural gas pipelines and offshore platforms. Finally, we compare the solutions and costs between the REORIENT and investment-only models. The investment-only model is the REORIENT model without the retrofit and abandonment planning functions. All data used in the test instances is provided at \cite{reorient_data}.

\subsection{Case study}
\begin{figure}[!thb]
    \centering
    \includegraphics[scale=0.9]{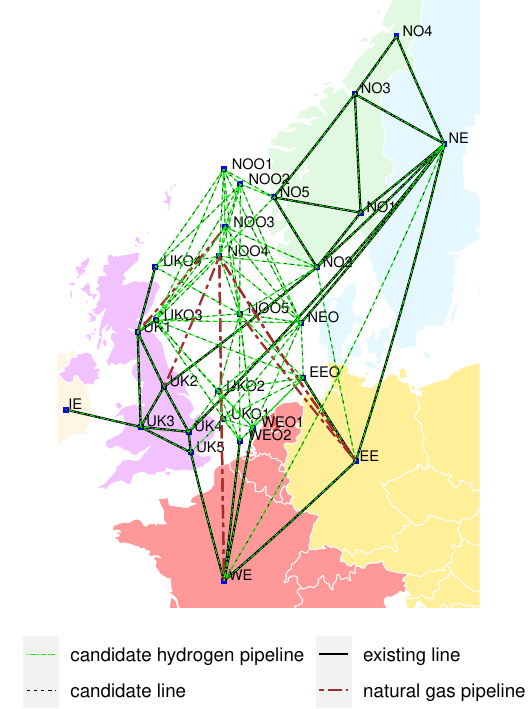}
    \caption{Illustration of the considered European energy system. The considered system includes 27 regions (each region can deploy 36 technologies), 87 transmission lines, 7 existing natural gas pipelines that can be retrofitted for hydrogen transport (some are overlapped), and 87 candidate new hydrogen pipelines.}
    \label{fig:ns_grid}
\end{figure}

We demonstrate the REORIENT model on the integrated strategic planning of the European energy system. The network topology is shown in Figure \ref{fig:ns_grid}. 
We make investment planning towards 2050 with a 5-year planning step. We implemented the algorithm and model in Julia 1.8.2 using JuMP \citep{jump} and solved with Gurobi 10.0 \citep{gurobi}. The problem instances contain up to 13 million continuous variables, 1860 binary variables, and 30 million constraints. We run the code on nodes of a computer cluster with a 2x 3.6GHz 8 core Intel Xeon Gold 6244 CPU and 384 GB of RAM, running on CentOS Linux 7.9.2009. 
% The Julia code and data for the case study have been made publicly available at \citep{reorient2022}. 
The parameters for the price process for oil and gas prices are presented in Table \ref{table:price_process_parameters}. 

\begin{table}[!htb]
\centering
\caption{ Estimated parameters for the price process, taken from \cite{Bakker2021MatureProgramming}.}
\label{table:price_process_parameters}
    \resizebox{8cm}{!}{%
    \begin{tabular}{ccccccc}
        \toprule
        & $\kappa$ & $\sigma_{\chi}$ & $\lambda_{\chi}$ & $\sigma_{\xi}$ &$\mu^{*}_{\xi}$ &$\rho_{\chi \xi}$\\ \hline
        Estimate    &0.407  &0.273  &-0.147 &0.149  &-0.007 &0.306\\
    \bottomrule
    \end{tabular}
    }
\end{table}

\begin{table}[!htb]
\centering
    \caption{Existing natural gas pipelines considered in the case study and their potential hydrogen transport capacity.}
    \label{table:Existing natural gas pipelines considered the case study.}
    \resizebox{10cm}{!}{%
        \begin{tabular}
            {llllS[table-format=2.0]}
            \toprule
            Model name  &   Name    & From      & To    & {Capacity (Ktonne/hour)}\\\hline
            Pipeline 1  &   Vesterled	&	NOO3	&	UK1	&	0.46\\
            Pipeline 2  &   Langeled	&	NOO4	&	UK2	&	0.98\\
            Pipeline 3  &   Zeepipe 1	&	NOO4	&	WE	&	0.58\\
            Pipeline 4  &   Franpipe	&	NOO4	&	WE	&	0.75\\
            Pipeline 5  &   Norpipe	    &	NOO5	&	EE	&	0.61\\
            Pipeline 6  &   Europipe 1	&	NOO4	&	EE	&	0.69\\
            Pipeline 7  &   Europipe 2	&	NOO4	&	EE	&	0.92\\
            \bottomrule
        \end{tabular}
    }
\end{table}

We use Gurobi as the base solver. We use the dual simplex algorithm to solve the RMP due to its relatively small size. The parameter \texttt{DegenMoves} has been turned on because we notice degeneracy makes the solver slow. We use the Barrier algorithm to solve the centred point problem to obtain a centred point. If \texttt{Presolve} is off and \texttt{Crossover} is off, then Gurobi will give a centred point. However, we turn on \texttt{Presolve} to reduce the problem size further. In addition, considering the scale of the problem, we choose to solve all the following instances to 1\% convergence tolerance.

\subsubsection{Computational results}
This section presents an overview of the problem instances and a performance analysis of the proposed algorithm. An overview of the test instances is presented in Table \ref{table:Overview of the cases used in the computational study.}. In the test instances, we consider operational problems with hourly time resolution. The test instances vary in the number of operational hours in each short-term scenario, short-term scenarios, and long-term scenarios. The problem instances have seven stages, which makes the problem instances large even with a few realisations of the parameters in each stage. The computational time is given in Tables \ref{table:Computational time of level set stabilised Benders} and \ref{table:Computational time of centred point stabilised Benders}, and note that Gurobi can only solve Case 1.

\begin{table}[!htb]
    \caption{Overview of the cases used in the computational study.}
    \label{table:Overview of the cases used in the computational study.}
        % \color{blue}
    \resizebox{\columnwidth}{!}{%
        \begin{tabular}
            {cS[table-format=3.0]S[table-format=1.0]S[table-format=2.0]|S[table-format=3.0]S[table-format=3.0]|S[table-format=1.1e1]S[table-format=4.0]S[table-format=1.1e1]|S[table-format=1.1e1]S[table-format=4.0]S[table-format=1.1e1]|S[table-format=1.1e1]S[table-format=1.1e1]}
            \toprule
            & {Operational periods} & \multicolumn{1}{c}{Short-term} & \multicolumn{1}{c}{Long-term} & \multicolumn{2}{c}{Number of decision nodes} & \multicolumn{3}{c}{Problem size (undecomposed)} & \multicolumn{3}{c}{Initial RMP size} & \multicolumn{2}{c}{Subproblem size}\\
            & {per short-term scenario} & \multicolumn{1}{c}{scenarios} & \multicolumn{1}{c}{scenarios} & \multicolumn{1}{c}{Operational nodes} & \multicolumn{1}{c}{Investment nodes} & \multicolumn{1}{c}{Continuous variables} & \multicolumn{1}{c}{Binary variables} & \multicolumn{1}{c}{Constraints} & \multicolumn{1}{c}{Continuous variables} & \multicolumn{1}{c}{Binary variables} & \multicolumn{1}{c}{Constraints} & \multicolumn{1}{c}{Continuous variables} & \multicolumn{1}{c}{Constraints} \\ \hline
Case	1	&	96	&	4	&	1	&	6	&	6	&	7.0e5	&	180	&	1.6e6	&2.0e4 &180 &2.3e4 &1.1e5 &1.7e5 \\
Case	2	&	672	&	4	&	1	&	6	&	6	&	4.8e6	&	180	&	1.1e7	&2.0e4 &180 &2.3e4 &8.0e5 &1.2e6\\
Case	3	&	96	&	4	&	53	&	114	&	62	&	1.3e7	&	1860	&	3.0e7	&3.0e5 &1860 &3.7e5 &1.1e5 &1.7e5\\
Case	4	&	96	&	4	&	49	&	107	&	58	&	1.2e7	&	1770	&	2.8e7	&2.9e5 &1770 &3.5e5 &1.1e7 &1.7e5\\
Case	5	&	96	&	8	&	49	&	107	&	58	&	8.5e7	&	1770	&	1.9e8	&2.9e5 &1770 &3.5e5 &2.3e5 &3.3e5\\
Case	6	&	96	&	4	&	95	&	273	&	179	&	3.2e7	&	5370	&7.1e7	&7.7e5 &5370 &9.4e5 &1.1e5 &1.7e5 \\
            \bottomrule
        \end{tabular}
    }
\end{table}

\begin{table}[!htb]
    \caption{Computational time of level set stabilised Benders. (Iters: iterations, Evals: subproblem evaluations)}
    \label{table:Computational time of level set stabilised Benders}
    \centering
    \resizebox{\columnwidth}{!}{%
    % \color{blue}
    \begin{tabular}{lrS[table-format=2.2]S[table-format=2.2]S[table-format=2.2]S[table-format=2.2]S[table-format=4.2]S[table-format=4.2]}
        \toprule
        & {Iters/Evals} &{Total time spent (h)} & {Master problem (\%)} & {Stabilisation problem (\%)} & {Subproblems and adaptive oracles (\%)} &{Lower bound (B€)} &{Upper bound (B€)}\\ \hline
Case 1	&	714/3615	&	5.53	&	16.34	&	43.06	&	40.58	&	1709.32	&	1726.29	\\
Case 2	&	894/4249	&	31.04	&	4.02	&	12.76	&	83.21	&	1882.78	&	1901.23	\\
        \bottomrule
        \multicolumn{8}{r}{This method did not solve Cases 3-6 after more than ten days of running.}
    \end{tabular}
    }
\end{table}

\begin{table}[!htb]
    \caption{Computational time of centred point stabilised Benders. (Iters: iterations, Evals: subproblem evaluations)}
    \label{table:Computational time of centred point stabilised Benders}
    \centering
    \resizebox{\columnwidth}{!}{%
    % \color{blue}
    \begin{tabular}{lrS[table-format=2.2e1]S[table-format=1.2e1]S[table-format=2.2e1]S[table-format=2.2e1]S[table-format=4.2]S[table-format=4.2]}
        \toprule
        & {Iters/Evals} & {Total time spent (h)} & {Master problem (\%)} & {Stabilisation problem (\%)} & {Subproblems and adaptive oracles (\%)} & {Lower bound (B€)} & {Upper bound (B€)}\\ \hline
Case 1	&	152/641	&	0.37	&	8.42	&	3.35	&	88.08	&	1717.62	&	1734.30	\\
Case 2	&	119/516	&	4.54	&	0.44	&	0.20	&	99.35	&	1887.87	&	1906.69	\\
Case 3	&	161/6502	&	49.61	&	9.26	&	8.88	&	81.85	&	1677.74	&	1694.62	\\
Case 4	&	159/5340	&	30.38	&	4.97	&	10.87	&	84.15	&	1378.83	&	1388.65	\\
Case 5	&	110/4371	&	55.12	&	0.87	&	2.87	&	96.25	&	1488.47	&	1501.71	\\
Case 6	&	134/9653	&	78.65	&	5.48	&	4.62	&	89.90	&	1329.80	&	1343.11	\\
        \bottomrule
    \end{tabular}
    }
\end{table}
By comparing Tables \ref{table:Computational time of level set stabilised Benders} and \ref{table:Computational time of centred point stabilised Benders}, we can see that by utilising the centred point, we reduce the computational time significantly. By comparing the percentage of the time spent on solving the stabilisation problem, we can see that solving CP takes much less of the total time than solving a quadratic programming stabilisation problem. We can also see that as we increase the number of strategic nodes, the percentage of time spent on the RMP and CP increases in both algorithms. This is because we add one cut per node per iteration, and as we have more nodes, the RMP and CP grow faster every iteration. Also, by comparing Cases 2 and 3, we observe no significant difference in the number of iterations, but the increase in time is more significant. This is because Cases 2 and 3 have the same amount of strategic decisions, but the subproblem in Case 3 is larger. In addition, there are significantly more operational nodes in Case 6 than the other cases which lead to much more subproblem evaluations. This eventually leads to a longer solution time.

\subsubsection{Sensitivity analysis}
\label{sec:sensitivity_analysis}
In this section, we use Case 3 to conduct a sensitivity analysis on the fixed retrofitting cost of pipelines and platforms. In addition, we also present the results of the investment decisions for a future energy system with a large amount of green and blue hydrogen production and transportation. 

We first conduct a sensitivity analysis on the retrofitting cost of natural gas pipelines. To this end, we consider two cases: Case A,  oil and gas production has stopped, and there is no natural gas transportation value in the pipelines, and Case B, oil and gas production is ongoing, and the natural gas pipelines are used for natural gas transport. Using Case A, our motivation is to understand under what cost range it would be more beneficial to retrofit natural gas pipelines that are not in operation compared with building new hydrogen pipelines. Case B is a more realistic case because most of the pipelines in the North Sea are in operation and have an export role. Using Case B, we want to analyse if retrofit to hydrogen will occur if that means the loss of oil and gas export profit. According to \cite{siemens_2021}, the cost of retrofitting the pipelines can be estimated at around 10-15\% of the new construction. Here, we conduct sensitivity on the cost from 5\%-30\% with a 5\% step. 

From Table \ref{table:Results of the expected retrofitting decisions in Case A.}, for pipelines 2-7, if the retrofitting cost is below 15\% of newly built cost, they will be retrofitted in all scenarios. However, for pipeline 1, when the retrofitting cost is less than 15\% of building a new one, it is retrofitted in 28 scenarios, and the retrofitting takes place in the third strategic stage. Pipeline 5 is only retrofitted at the end of the planning horizon when the cost is 15\% of building a new one. When the retrofitting cost is higher than 20\%, some of the pipelines are not retrofitted. Instead, the model decides to build new pipelines and a different network topology to achieve the minimum cost. We can see that different oil and gas price scenarios affect the retrofitting decisions. 
\begin{table}[!htb]
\centering
    \caption{Results of the expected retrofitting decisions in Case A.}
    \label{table:Results of the expected retrofitting decisions in Case A.}
    \resizebox{16cm}{!}{%
        \begin{tabular}
            {lrrrrrrr}
            \toprule
            Cost (\% of new one)	&	{Pipeline 1}	&	{Pipeline 2}	&	{Pipeline 3}	&	{Pipeline 4}	&	{Pipeline 5}	&	{Pipeline 6}	&	{Pipeline 7}	\\ \hline
            5\%	&	(0, 0, 2, 4, 8, 16, 28)	&	$*$	&	$*$	&	$*$	&	$*$	&	$*$	&	$*$	\\
            10\%	&	(0, 0, 2, 4, 8, 16, 28)	&	$*$	&	$*$	&	$*$	&	$*$	&	$*$	&	$*$	\\
            15\%	&	(0, 0, 2, 4, 8, 16, 28)	&	$*$	&	$*$	&	$*$	&	(0, 0, 0, 0, 0, 2, 3)	&	$*$	&	$*$	\\
            20\%	&	$-$	&	$*$	&	$*$	&	$-$	&	$-$	&	$*$	&	$*$	\\
            25\%	&	$-$	&	$-$	&	$*$	&	$*$	&	(0, 0, 0, 2, 4, 7, 11)	&	$*$	&	$*$	\\
            30\%	&	$-$	&	$-$	&	$*$	&	$-$	&	(0, 0, 0, 0, 2, 4, 7)	&	$*$	&	$*$	\\
            \bottomrule
            \multicolumn{8}{r}{$*$: the pipeline is retrofitted in all strategic nodes, $-$: the pipeline is not retrofitted.}\\
            \multicolumn{8}{r}{$(\{x_i, i=1,..,7\}):$ the number of decision nodes, $x_i$, that retrofitting of the pipeline takes place in stage $i$.}
        \end{tabular}
    }
\end{table} 

\begin{table}[!htb]
\centering
    \caption{Results of the expected retrofitting decisions in Case B.}
    \label{table:Results of the expected retrofitting decisions in Case B.}
    \resizebox{\columnwidth}{!}{%
        \begin{tabular}
            {lrrrrrrr}
            \toprule
            Cost (\% of new one)	&	{Pipeline 1}	&	{Pipeline 2}	&	{Pipeline 3}	&	{Pipeline 4}	&	{Pipeline 5}	&	{Pipeline 6}	&	{Pipeline 7}	\\ \hline
            5\%	&	$-$	&	(0, 0, 0, 0, 2, 4, 6)	&	(0, 0, 0, 2, 4, 8, 14)	&	(0, 0, 0, 0, 2, 4, 7)	&	(0, 0, 0, 0, 2, 4, 8)	&	(0, 0, 0, 0, 2, 4, 7)	&	(0, 0, 0, 0, 2, 4, 7)	\\
            10\%	&	$-$	&	$-$	&	(0, 0, 0, 2, 4, 8, 14)	&	(0, 0, 2, 4, 8, 16, 28)	&	(0, 0, 0, 0, 0, 2, 4) 	&	(0, 0, 0, 2, 4, 8, 14)	&	(0, 0, 0, 2, 4, 8, 14)	\\
            15\%	&	$-$	&	$-$	&	(0, 0, 0, 2, 4, 8, 14)	&	(0, 0, 0, 2, 4, 8, 14)	&	(0, 0, 0, 0, 2, 4, 6)	&	(0, 0, 0, 2, 4, 8, 14)	&	(0, 0, 0, 0, 2, 4, 7)	\\
            20\%	&	$-$	&	$-$	&	(0 ,0, 0, 0, 2, 4, 8)	&	(0, 0, 0, 0, 2, 4, 6)	&	(0, 0, 0, 0, 2, 4, 6)	&	(0, 0, 0, 0, 2, 4, 8)	&	(0, 0, 0, 0, 2, 4, 8)	\\
            25\%	&	$-$	&	$-$	&	(0, 0, 0, 2, 4, 8, 14)	&	(0, 0, 0, 0, 2, 4, 6)	&	(0, 0, 0, 0, 2, 4, 6)	&	(0, 0, 0, 2, 4, 8, 14) 	&	(0, 0, 0, 2, 4, 8, 14)	\\
            30\%	&	$-$	&	$-$	&	$-$	&	$-$	&	$-$	&	(0, 0, 0, 0, 2, 4, 8)	&	(0, 0, 0, 0, 2, 4, 6)	\\
            \bottomrule
\multicolumn{8}{r}{$-$: the pipeline is not retrofitted.}\\
\multicolumn{8}{r}{$(\{x_i, i=1,..,7\}):$ the number of decision nodes, $x_i$, that retrofitting of the pipeline takes place in stage $i$.}
        \end{tabular}
    }
\end{table} 

From Table \ref{table:Results of the expected retrofitting decisions in Case B.}, it can be observed, compared with Case A, that the economic viability of pipeline retrofit is harder if the pipelines are already used for natural gas transport. However, most pipelines are still retrofitted for hydrogen transportation in later investment stages. From Tables \ref{table:Results of the expected retrofitting decisions in Case A.} and \ref{table:Results of the expected retrofitting decisions in Case B.}, we can see that retrofit decisions are sensitive to the retrofitting cost, and oil and gas prices. Also, retrofit sometimes only take place in specific price scenarios.  

Secondly, we conduct a sensitivity analysis on the retrofitting cost of oil and gas platforms. By doing so, we aim to analyse: (1) if retrofitting can help delay or even avoid the costly abandonment campaign and (2) understand the relation between retrofitting existing platforms for OEHs and building new OEHs. We assume that the fixed part of the retrofitting cost is half of the removal cost, and conduct sensitivity around this cost. For each platform cluster, we consider a fixed part of the retrofitting cost ranging from €10 million to €2 billion. However, the results suggest that it is not economical to retrofit platforms for hydrogen production under this price range due to the massive loss of oil and gas export profit. The model decides to conduct an abandonment campaign for all platform clusters by the end of the planning horizon. This means that based on the cost models that are used, retrofitting platforms for hydrogen production is more costly than abandonment. Also, due to the oil and gas export profit, the platforms will produce as long as possible until they must be abandoned. In this case study, the platforms must be retrofitted or abandoned by 2050. This suggests that repurposing platforms for other use may need stronger incentives in addition to economic factors. 

\begin{figure}[!thb]
    \centering
    \includegraphics[scale=0.8]{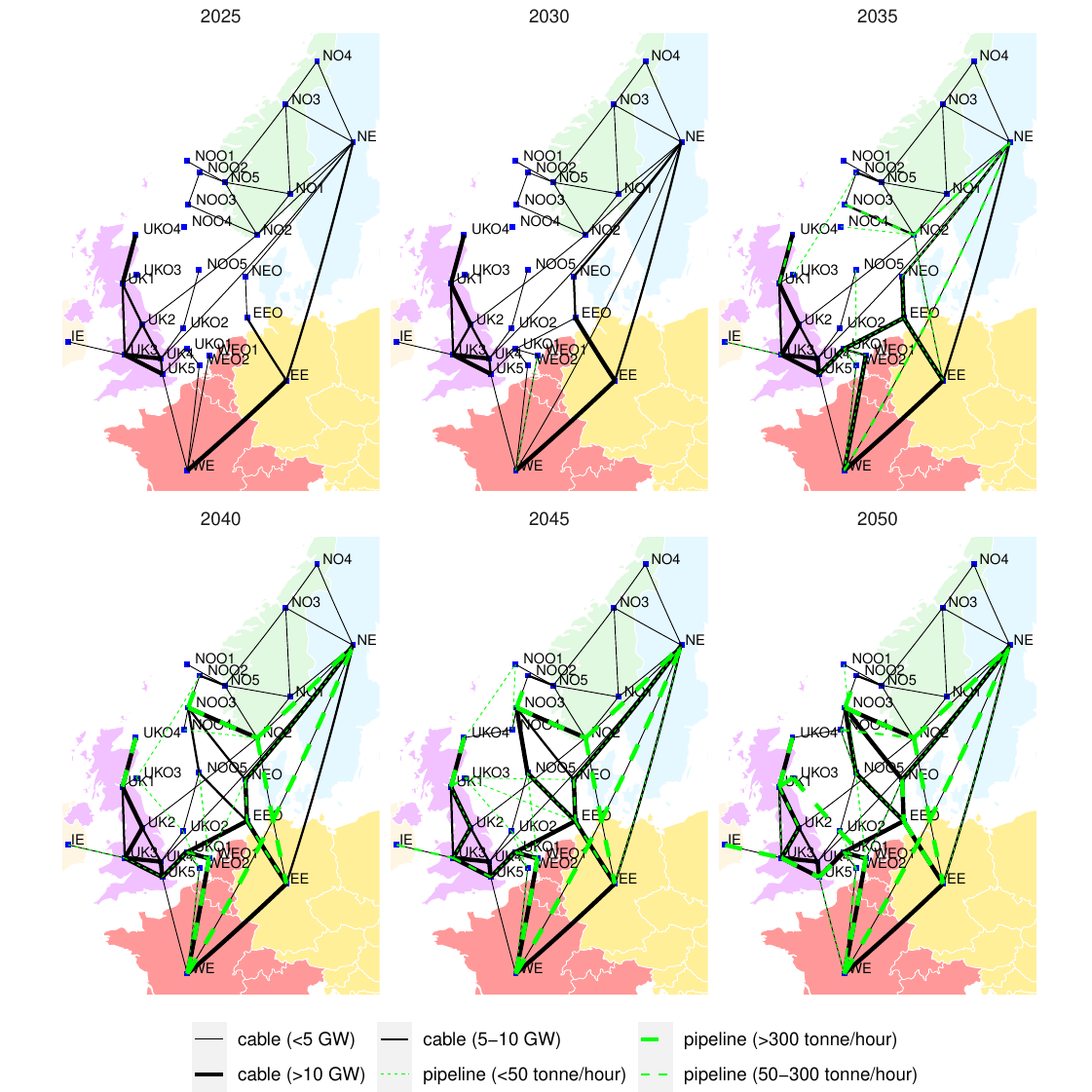}
    \caption{Expected solution of the grid design towards 2050 (investment only model).}
    \label{fig:grid_investment_only}
\end{figure}
\subsubsection{Comparison between the REORIENT model and an investment planning only model}
In this section, we use Case 3 and analyse the difference between an investment-planning-only model and the proposed integrated model regarding investment decisions and costs. We fix the retrofitting cost of pipelines to 15\% of the cost of its newly built counterpart. In the following, we report the results of expected strategic decisions regarding the grid design and capacities of the technologies in each decision stage. 

\begin{figure}[!thb]
    \centering
    \includegraphics[scale=0.8]{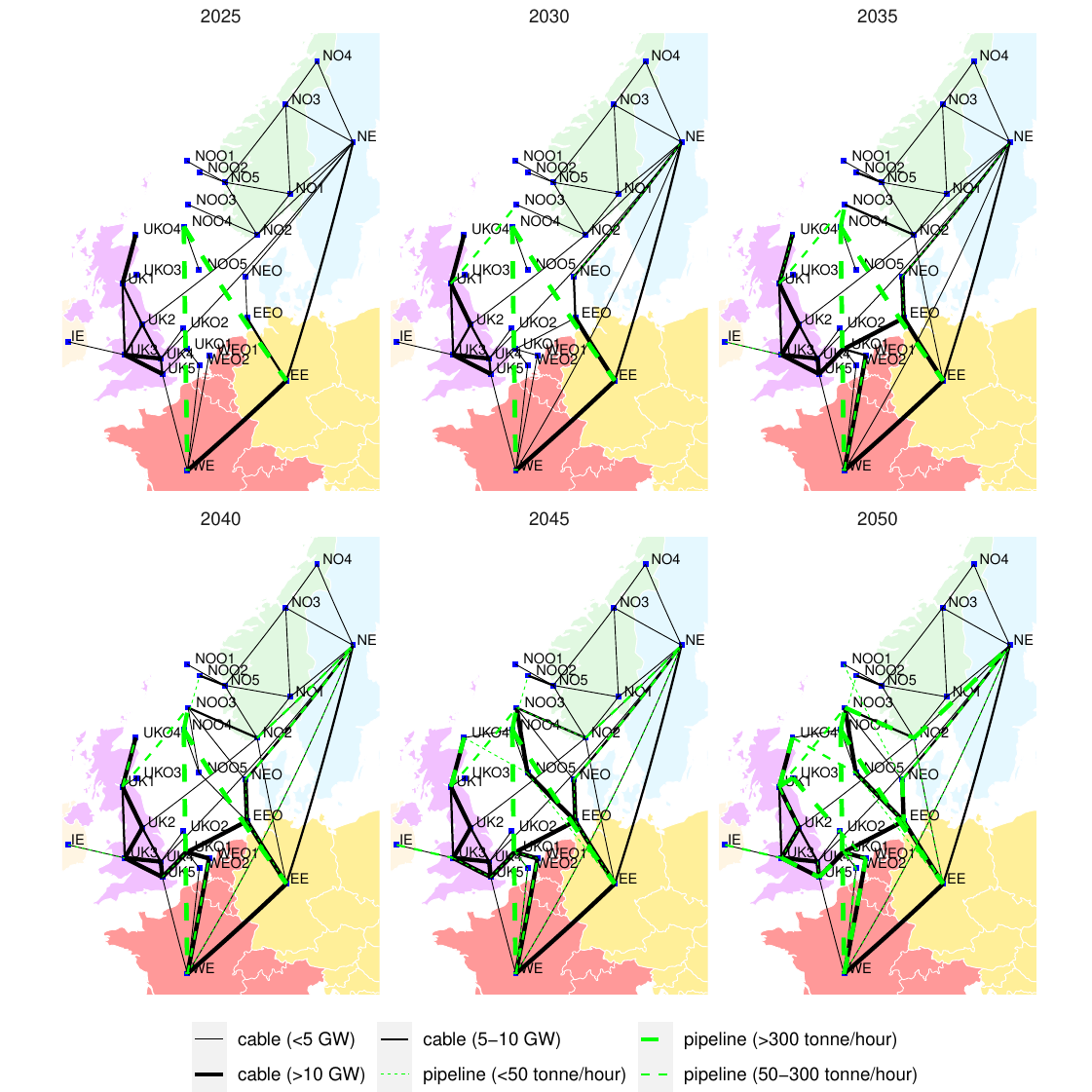}
    \caption{Expected solution of the grid design towards 2050 (REORIENT model).}
    \label{fig:grid_REORIENT}
\end{figure}

\begin{figure}[!thb]
    \centering
    \includegraphics[scale=0.9]{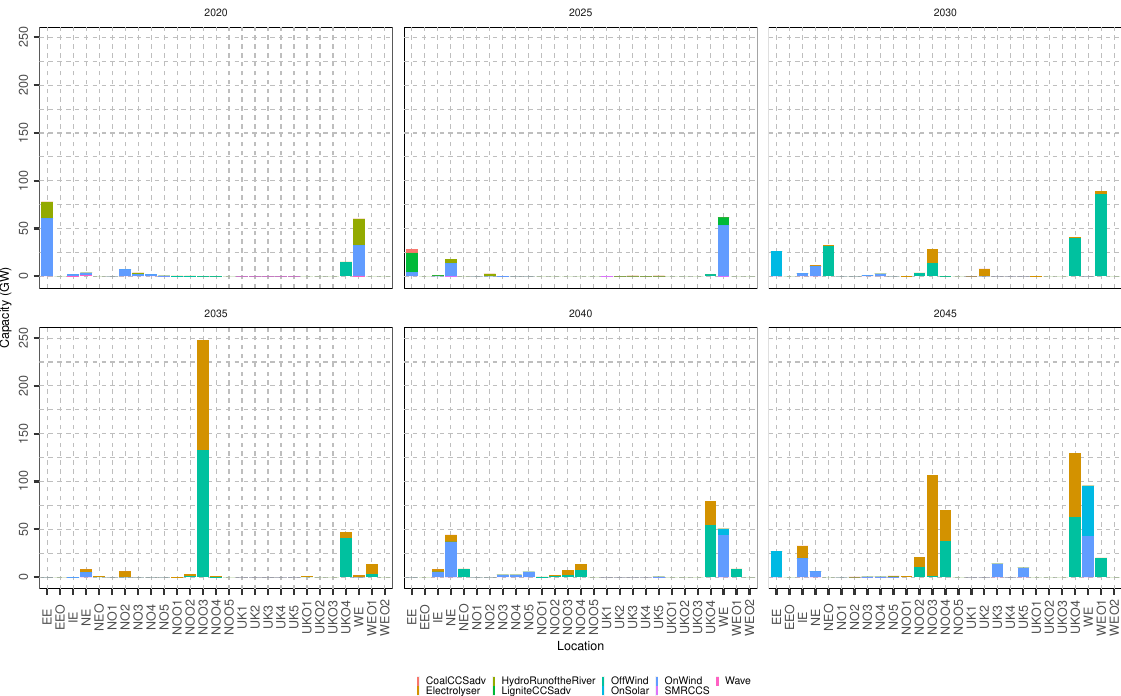}
    \caption{Expected investment decisions towards 2050 (investment only model).}
    \label{fig:tech_investment_only}
\end{figure}
\begin{figure}[!thb]
    \centering
    \includegraphics[scale=0.9]{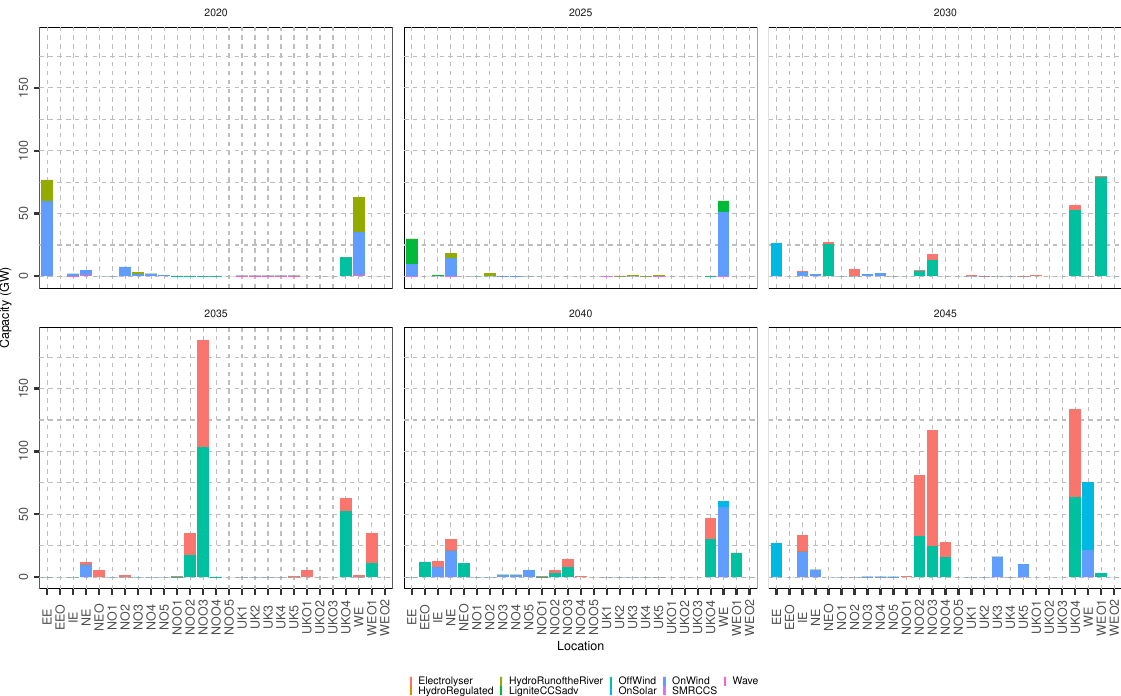}
    \caption{Expected investment decisions towards 2050 (REORIENT model).}
    \label{fig:tech_REORIENT}
\end{figure}

From Figures \ref{fig:grid_investment_only} and \ref{fig:grid_REORIENT}, we can see that the network topology is noticeably different. By 2050, there will be 32 pipelines built compared with 28 pipelines in the investment-only model. The line connecting NE and NO1 has less capacity in the REORIENT model compared with the investment-only model. In both cases, the UK onshore power system transmission is reinforced, however, a 3 GW difference in the line connecting UK3 and UK5 is observed, followed by a 2 GW difference in UK4-UK5. By 2050, the line NEO-EEO will have 44.90 GW capacity in REORIENT model compared with 36.81 GW in its counterpart. NOO3-NEO presents a significant difference as well with 4.4 GW capacity in REORIENT model and 12.26 GW in the investment-only model by 2050. NOO2 and NOO3 are not connected in the REORIENT model but are connected in the other model. 

By comparing Figures \ref{fig:tech_investment_only} and \ref{fig:tech_REORIENT}, we notice that in both models,  NOO3 is an important offshore region which receives significant investment in offshore wind and electrolysers due to its location and high wind availability. A major difference is found in offshore wind capacity in NEO and UKO4. The REORIENT model has a higher investment in offshore wind in NEO in all investment steps. 

\begin{table}[!htb]
    \centering
    \caption{Results of accumulated capacity in Europe (investment-only model).}
    \label{table:Results of accumulated capacity in Europe (investment-only model).}
    \resizebox{\columnwidth}{!}{%
    \begin{tabular}{lS[table-format=3.2]S[table-format=3.2]S[table-format=3.2]S[table-format=3.2]S[table-format=3.2]S[table-format=4.2]S[table-format=2.2]}
    \toprule
     {Year} & {Offshore wind}	&	{Onshore wind}	&	{Onshore solar} &		{Electrolyser offshore} &		{Electrolyser onshore}	&	{Transmission line}	&	{Hydrogen pipeline}\\
    & {(GW)}	&	{(GW)}	&	{(GW)} &		{(GW)} &		{(GW)}	&	{(GW)}	&	{(ktonne/hour)}\\ \hline
     2025	&	35.72	&	279.24	&	119.83	&	0.00	&	0.00	&	148.66	&	0.00	\\
    2030	&	36.26	&	354.17	&	119.83	&	0.00	&	0.00	&	190.46	&	0.00	\\
    2035	&	212.08	&	364.02	&	146.29	&	12.46	&	8.66	&	586.61	&	0.73	\\
    2040	&	397.01	&	373.98	&	146.29	&	161.64	&	14.74	&	909.83	&	6.98	\\
    2045	&	465.95	&	355.11	&	150.86	&	174.57	&	20.51	&	985.28	&	7.80	\\
    2050	&	608.14	&	269.46	&	171.42	&	248.01	&	27.19	&	1084.87	&	15.70	\\
    \bottomrule
    \end{tabular}
    }
\end{table}

\begin{table}[!htb]
    \centering
    \caption{Results of accumulated capacity in Europe (REORIENT model).}
    \label{table:Results of accumulated capacity in Europe (REORIENT model).}
    \resizebox{\columnwidth}{!}{%
    \begin{tabular}{lS[table-format=3.2]S[table-format=3.2]S[table-format=3.2]S[table-format=3.2]S[table-format=3.2]S[table-format=4.2]S[table-format=2.2]}
    \toprule
     {Year} & {Offshore wind}	&	{Onshore wind}	&	{Onshore solar} &		{Electrolyser offshore} &		{Electrolyser onshore}	&	{Transmission line}	&	{Hydrogen pipeline}\\
    & {(GW)}	&	{(GW)}	&	{(GW)} &		{(GW)} &		{(GW)}	&	{(GW)}	&	{(ktonne/hour)}\\ \hline
    2025	&	35.92	&	277.39	&	119.83	&	0.00	&	0.00	&	148.84	&	2.94	\\
    2030	&	38.33	&	348.67	&	119.83	&	0.00	&	0.00	&	190.73	&	3.19	\\
    2035	&	215.35	&	366.35	&	146.29	&	20.25	&	9.35	&	536.56	&	3.59	\\
    2040	&	394.34	&	371.60	&	146.29	&	156.90	&	20.38	&	904.63	&	6.19	\\
    2045	&	461.80	&	357.82	&	151.91	&	172.95	&	21.44	&	998.00	&	6.93	\\
    2050	&	593.89	&	293.15	&	171.42	&	253.37	&	24.85	&	1061.02	&	13.63	\\
    \bottomrule
    \end{tabular}
    }
\end{table}

Tables \ref{table:Results of accumulated capacity in Europe (investment-only model).} and \ref{table:Results of accumulated capacity in Europe (REORIENT model).} present the accumulated capacity of each technology in each region offshore wind will surpass onshore wind and become the most important renewable power supply by 2050. In the investment-only model, more than 23 GW more transmission line capacity is needed by 2050 compared with the results using the REORIENT model. In both models, hydrogen is produced mainly from SMRCCS at the initial stages but gradually replaced by electrolysers. Also, both models decide to produce green hydrogen mainly offshore. The hydrogen pipeline capacity is lower in the REORIENT model compared with the counterpart by the end of the planning stage. 

In addition, the total cost over the given planning horizon is €1694.47 billion in the investment-only model and €1691.50 billion in the REORIENT model. Furthermore, the REORIENT yields 24\% lower investment cost in the North Sea region compared with the traditional investment-only model. This shows the potential value of doing integrated planning. The value of REORIENT model can be further revealed once more retrofitting options are included, e.g., by including more existing natural gas pipelines. 

\section{Discussion}
\label{sec:discussion}
In this paper, we integrated investment, retrofit and abandonment planning in a multi-horizon stochastic MILP model. The model is generally applicable to studying a specific planning problem for a production plant or a large-scale energy system planning problem for a region.

We used the model to study the investment planning of a European energy system. We considered regional retrofit at a high level and conduct a techno-economical analysis. Unlike traditional retrofit models for process systems, we have omitted detailed modelling of the processes, which is a compromise due to the large scale of the study. The sensitivity analysis presented in this paper can be used as a benchmark for future studies. 

In the case study, we find that although reducing retrofitting costs can trigger the retrofitting of some oil and gas infrastructures, it may not be a sufficiently strong incentive for platform retrofitting compared with pipeline retrofitting. This is because the loss of oil and gas profit is much larger than the reduction in retrofitting cost. Additionally, for platforms, a lot of investment needs to be made for producing green hydrogen upon removing the existing structure. Other driving factors, such as policies, are therefore needed to encourage oil and gas operators to retrofit their infrastructure or reduce production for the energy transition. 

From a computational perspective, the proposed algorithm solved the problem instances efficiently. The problem instances have many regions and technologies and, therefore, are highly degenerate. The CP helps the proposed stabilised adaptive Benders algorithm to converge faster. In addition, the proposed stabilised adaptive Benders decomposition can be applied to a class of problems that can be formulated as Equations \eqref{eq:MHSP} and \eqref{eq:benders_subproblem_general}. Other strategies to accelerate Benders decomposition, including adding combinatorial cuts, trust region, local branching methods, and partial surrogate are tested. However, the improvement in performance is not significant. 

There are some limitations of the case study: (1) in the case study, the offshore fields are aggregated into representative fields, which loses the modelling of the retrofit and abandonment for specific fields; (2) there are other parameters that may affect the investment decisions such as CO$_2$ budget can be an uncertain parameter and change the results. We believe that oil and gas prices have a more direct relation and economic trade-off with retrofit and abandonment decisions; therefore, we choose to consider oil and gas prices as the uncertain parameter, and (3) we only consider green and blue hydrogen. However, hydrogen produced by other means may also be relevant and affect the results. In addition, although the MHSP can include uncertainty from short-term and long-term time horizons more efficiently, uncertain parameters such as oil and gas prices can be affected by regional decisions taken. This can not be captured using our MHSP model. However, North Sea oil and gas prices are mainly affected by global actions not regional actions so our model is a good approximation. If local decisions do have a significant effect on oil and gas prices, there have been modelling and computational strategies for multi-stage stochastic programming with endogenous and exogenous uncertainties \citep{Goel2006AUncertainty, Apap2017ModelsUncertainties} and it may be possible to combine these approaches with MHSP to address this limitation.

\section{Conclusions and future work}
\label{sec:conclusions}
This paper has presented the REORIENT model, a multi-horizon stochastic MILP for integrated investment, retrofit and abandonment energy-system planning. The major novelties and contributions are: (1) we developed an MHSP model for integrated investment, retrofit and abandonment planning of energy systems, (2) we included uncertainty from both strategic and operational time horizons in such a model, (3) a centred point stabilised adaptive Benders decomposition algorithm was developed to solve large-scale MILP faster, (4) we extended the centred point stabilisation, which was used for linear programs, to solve MILP problems and provided convergence proof, and (5) the triggering parameters for retrofitting is investigated by conducting sensitivity analysis and a comparison between the REORIENT model and investment planning only model is made. Results from our case study indicate that: (1) for pipelines that are not in use, when the retrofitting cost is below 20\% of the cost of building new ones, it is more economical to retrofit most of the pipelines than building new ones. For pipelines that transport natural gas, it is economical to be retrofitted in some natural gas price scenarios, (2) platform clusters keep producing oil and gas rather than being retrofitted for hydrogen use, and the clusters abandonment takes place at the last investment stage, (3) compared with an investment planning model, the REORIENT model yields €3 billion lower total cost, and 24\% lower investment cost in the North Sea region, and (4) the proposed Benders algorithm can solve the model efficiently. For the smaller cases, it is 6.8 times faster than the level method stabilised Benders decomposition, and it can solve all the larger cases which the level method stabilised Benders which could not solve.

In the future, the REORIENT model can be used for more energy systems analysis, such as investigating the integrated planning for other regions, such as the continental shelf of the United States, or focusing on some specific platforms in a smaller region. In addition, decomposition methods, such as Lagrangean type algorithms and progressive hedging algorithms, could be applied to solve the reduced master problems should this becomes a bottleneck in solving larger stochastic instances.

\section*{CRediT author statement}
\textbf{Hongyu Zhang:} Conceptualisation, Methodology, Software, Validation, Formal analysis, Investigation, Visualisation, Data curation, Writing - original draft, Writing - review \& editing. \textbf{Ignacio E. Grossmann:} Conceptualisation, Methodology, Supervision, Writing - review \& editing. \textbf{Ken McKinnon:} Conceptualisation, Methodology, Writing - review \& editing. \textbf{Rodrigo Garcia Nava:} Conceptualisation, Methodology, Software. \textbf{Brage Rugstad Knudsen:} Conceptualisation, Methodology, Supervision, Writing - review \& editing, Funding acquisition. \textbf{Asgeir Tomasgard:} Conceptualisation, Methodology, Supervision, Writing - review \& editing, Funding acquisition.

\section*{Declaration of competing interest}
The authors declare that they have no known competing financial interests or personal relationships that could have appeared to influence the work reported in this paper.

\section*{Acknowledgements}
This work was supported by the Research Council of Norway through PETROSENTER LowEmission [project code 296207].




 % \footnote{For further information, please refer to \url{https://docs.julialang.org/en/v1/manual/asynchronous-programming/}}

% \section*{References}
%\bibliographystyle{elsarticle-harv}
\clearpage

\setlength{\bibsep}{0pt plus 0.3ex}
\footnotesize{
\bibliographystyle{model5-names}
\bibliography{retrofit_offshore_system_3rd_r2_clean}}
\clearpage

\appendix

\section{Nomenclature}
\label{sec:nonmenclature}
\begin{description}[itemsep=-5pt, leftmargin=!,labelwidth=\widthof{\bfseries $X^{MaxReT/MinReT}_{pi}$}]
\item [\textbf{Investment planning model indices and sets}]
\item [$p \in \mathcal{P}$] set of technologies
\item [$p \in \mathcal{P}^{R}$] set of candidate retrofit technologies
\item [$p \in \mathcal{P}^{R}_{p}$] set of candidate technologies that an existing technology $p$ can be retrofitted to $(p \in \mathcal{P}^{R})$ which including abandonment and prolong
\item [$p \in \mathcal{P}^{RT}$] set of candidate technologies be retrofitted to
\item [$i \in \mathcal{I}^{Ope}$] set of operational nodes
\item [$i \in \mathcal{I}^{Inv}$] set of investment nodes
\item [$j \in \mathcal{I}^{Inv}_{i}$] set of investment nodes $j$ $(j \in \mathcal{I}^{Inv})$ succeed to investment node $i$ $(i \in \mathcal{I}^{Inv})$
\item [$j \in \mathcal{I}^{Ope}_{i}$] set of operational nodes $j$ $(j \in \mathcal{I}^{Ope})$ succeed to investment node $i$ $(i \in \mathcal{I}^{Inv})$
\item [$(x, \theta, \lambda) \in \mathcal{F}_{i(k-1)}$] set of the Benders cut built up to iteration $k-1$, where $x$ is the vector of sampled points, $theta$ and $\lambda$ are the actually cost of subproblem at the sampled points, and the vector of subgradients at the sampled points, respectively. 

\item [\textbf{Investment planning model parameters}]
\item [$C^{InvV}_{pi}$] unitary investment cost of technology $p$ in investment node $i$ ($p \in \mathcal{P}, i \in \mathcal{I}^{Inv}$) [€/MW, €/MWh, €/kg]
\item [$C^{InvF}_{pi}$] fixed capacity independent investment cost of technology $p$ in investment node $i$ ($p \in \mathcal{P}, i \in \mathcal{I}^{Inv}$) [€]
\item [$C^{Fix}_{pi}$] unitary fix operational and maintenance cost of technology $p$ in operational node $i$ ($p \in \mathcal{P},i \in \mathcal{I}^{Ope}$) [€/MW, €/MWh, €/kg]
\item [$C^{ReTV}_{pi}$] unitary investment cost of retrofitted technology $p$ in investment node $i$ ($p \in \mathcal{P}^{RT}, i \in \mathcal{I}^{Inv}$) [€/MW, €/MWh, €/kg]
\item [$C^{ReTF}_{pi}$] fixed capacity independent investment cost of retrofitted to technology $p$ in investment node $i$ ($p \in \mathcal{P}^{RT}, i \in \mathcal{I}^{Inv}$) [€]
\item [$C^{ReTFixO}$] fix operational cost of the technology that is retrofitted to $p$ in investment node $i$ ($p \in \mathcal{P}^{RT}, i \in \mathcal{I}^{Inv}$) [€]
\item [$C^{ReFFixO}_{pi}$] fix operational cost of retrofitted technology $p$ in investment node $i$ ($p \in \mathcal{P}^{RT}, i \in \mathcal{I}^{Inv}$) [€]
\item [$X^{MaxInv/MinInv}_{pi}$] maximum/minimum built capacity of technology $p$ in investment node $i$ ($p \in \mathcal{P}, i \in \mathcal{I}^{Inv}$) [MW, MWh, kg]
\item [$X^{MaxAcc}_{p}$] maximum installed capacity of technology over the planning horizon $p$ ($p \in \mathcal{P}$) [MW, MWh, kg]
\item [$\kappa$] scaling effect depending on time step between successive investment nodes
\item [$H_p$] lifetime of technology $p$ ($p \in \mathcal{P}$)
\item [$X^{HistReF}_{pi}$] historical capacity of existing technology that can be retrofitted [MW, MWh, kg]
\item [$X^{Hist}_{pi}$] historical capacity of technology $p$ in operational node $i$ ($p \in \mathcal{P}, i \in \mathcal{I}^{Ope}$) [MW, MWh, kg]
\item [$X^{MaxReT/MinReT}_{pi}$] maximum/minimum built capacity of technology $p$ in investment node $i$ ($p \in \mathcal{P}^{RT}, i \in \mathcal{I}^{Inv}$) [MW, MWh, kg]
\item [$X^{MaxAccReT}_{pi}$] maximum installed capacity of technology $p$ ($p \in \mathcal{P}^{RT}$) [MW, MWh, kg]
\item [$x_{i}$] right hand side of the operational problem
\item [$c_{i}$] cost coefficients of the operational problem
\item [$\pi^{Inv/Ope}_{i}$] discount factor multiplied probability of investment/operational node $i$, ($i \in \mathcal{I}^{Inv}$/$i \in \mathcal{I}^{Ope}$)
\item [$\mu^{E}_{i}$] CO$_2$ budget at operational node $i$ ($i \in \mathcal{I}^{Ope}$)
\item [$\mu^{DP}_{i}$] scaling factor on power demand at operational node $i$ ($i \in \mathcal{I}^{Ope}$)
\item [$\mu^{P}_{i}$] scaling factor on oil and gas production at operational node $i$ ($i \in \mathcal{I}^{Ope}$)
\item [$\mu^{DHy}_{i}$] scaling factor on hydrogen demand at operational node $i$ ($i \in \mathcal{I}^{Ope}$)
\item [$C^{CO2}_{i}$] CO$_2$ emission price at operational node $i$ ($i \in \mathcal{I}^{Ope}$)
\item [$\mathcal{S}^{Ope}_{i}$] strategic stage of operational node $i$ $(i \in \mathcal{I}^{Ope})$
\item [$\mathcal{S}^{Inv}_{i}$] strategic stage of investment node $i$ $(i \in \mathcal{I}^{Inv})$

\item [\textbf{Investment planning model variables}]
\item [$x_{pi}^{Acc}$] accumulated capacity of device $p$ in operational node $i$ ($p \in \mathcal{P}, i \in \mathcal{I}^{Ope}$) [MW, MWh, kg]
\item [$x_{pi}^{Inv}$] newly invested capacity of device $p$ in investment node $i_{0}$ ($p \in \mathcal{P}, i \in \mathcal{I}^{Inv}$) [MW, MWh, kg]
\item [$y^{Inv}_{pi}$] 1 if technology $p$ is newly invested in investment node $i$, 0 otherwise ($p \in \mathcal{P}, i \in \mathcal{I}^{Inv}$)
\item [$y^{ReT}_{pi}$] 1 if technology $p$ is retrofitted to in investment node $i$, 0 otherwise ($p \in \mathcal{P}^{RT}, i \in \mathcal{I}^{Inv}$)
\item [$x^{AccReT}_{pi}$] accumulated capacity of technology $p$ that is retrofitted to in operational node $i$ ($p \in \mathcal{P}^{RT}, i \in \mathcal{I}^{Ope}$)
\item [$y^{ReF}_{pi}$] 1 if technology $p$ is retrofitted from in investment node $i$, 0 otherwise ($p \in \mathcal{P}^{R}, i \in \mathcal{I}^{Inv}$)
\item [$x^{AccReF}_{pi}$] accumulated capacity of retrofitted from technology in operational node $i$ ($p \in \mathcal{P}^{R}, i \in \mathcal{I}^{Ope}$)
\item [$x^{ReT}_{pi}$]  in operational node $i$ ($p \in \mathcal{P}^{R}, i \in \mathcal{I}^{Ope}$)
\item [$c^{INV}$] total investment and fixed operating and maintenance costs [€]
\item [$c^{OPE}_i$] approximated operational cost in operational node $i$ in Benders decomposition ($i \in \mathcal{I}^{Ope}$) [€]

\item [\textbf{Operational model indices and sets}]
\item [$n \in \mathcal{N}$] set of time slices
\item [$t \in \mathcal{T}$] set of hours in all time slices
\item [$t \in \mathcal{T}_{n}$] set of hours in time slice $n$ $(n \in \mathcal{N})$
\item [$l \in \mathcal{L}$] set of transmission lines
\item [$l \in \mathcal{L}^{Hy}$] set of hydrogen pipelines
\item [$l \in \mathcal{L}^{Out/In}_z$] set of transmission lines go out of/into region $z$
\item [$l \in \mathcal{L}^{HyOut/HyIn}_z$] set of hydrogen pipelines go out of/into region $z$
\item [$g \in \mathcal{G}$] set of thermal generation
\item [$r \in \mathcal{G}^{R}$] set of renewable generation
\item [$g \in \mathcal{G}^{H}$] set of hydropower generation including run of the river $\mathcal{G}^{HRor}$ and seasonal $\mathcal{G}^{HSea}$
\item [$s \in \mathcal{S}^{E}$] set of electricity storage
\item [$s \in \mathcal{S}^{Hy}$] set of hydrogen storage
\item [$b \in \mathcal{B}^{E}$] set of electric boilers
\item [$r \in \mathcal{R}$] set of SMRCCS
\item [$e \in \mathcal{E}$] set of electrolysers
\item [$f \in \mathcal{F}$] set of fuel cells
\item [$z \in \mathcal{Z}^{P}$] set of all platform clusters
\item [$z \in \mathcal{Z}$] set of all locations
\item [$p \in \mathcal{P}^{*}$] set of all thermal generators, electric boilers, electrolysers, electricity storage, fuel cells and seasonal hydropower generation $(\mathcal{P}^{*}=\mathcal{G} \cup \mathcal{B}^{E} \cup \mathcal{E} \cup \mathcal{S}^{E} \cup \mathcal{F} \cup \mathcal{G}^{HSea})$
\item [$v \in \mathcal{V}^{*}$] set of hydrogen storage and SMRCCS plants $(\mathcal{V}^{*}=\mathcal{S}^{Hy} \cup \mathcal{R})$ 

\item [\textbf{Operational model parameters}]
\item [$\mu^{E}$] CO$_2$ emission limit (tonne)
\item [$\mu^{DP/DH/DHy}$] scaling effect on power demand/heat demand/hydrogen demand
\item [$H_t$] number of hour(s) in one operational period $t$
\item [$\pi_t$] weighted length of one operational period $t$
\item [$R_{rt}^{GR}$] capacity factor of renewable unit $r$ in period $t$ ($r \in \mathcal{R}, t \in \mathcal{T}$)
\item [$\eta^*$] efficiency of electric boilers, fuel cells, thermal generators, electric storage and transmission lines $*=\{\text{BE, SE, L, HrG\}}$ indexed by related sets
\item [$E_{g}^{G}$]  CO$_2$ emission factor of thermal generation $g$ ($g \in \mathcal{G}$) [t/MWh]
\item [$C_g^G$] total operational cost of generating 1 MW power from thermal generation $g$ ($g \in \mathcal{G}$) [€/MW]
\item [$C^{Shed,l}$] load shed penalty cost of power $(l=P)$, heat $(l=H)$ and hydrogen $l=Hy$ [€/MW, €/kg]
\item [$\sigma^{Res}_{z}$] spinning reserve factor in region $z$ ($z \in \mathcal{Z}$)
\item [$\alpha^{G}_{g}$] maximum ramp rate of generators ($g \in \mathcal{G}$) [MW/MW]
\item [$\rho^{E}$] conversion factor of electrolyser to hydrogen [MWh/kg]
\item [$P^{DP/DH}_{zt}$] power demand/heat demand in location $z$ period $t$ $(z \in \mathcal{Z}, t \in \mathcal{T})$ [MW]
\item [$\rho^{F}$] hydrogen consumption factor of fuel cell [kg/MW]
\item [$P_{g}^{AccG}$] accumulated capacity of thermal generator $g$ $(g \in \mathcal{G})$ [MW]
\item [$P_{g}^{AccHRor}$] accumulated capacity of run of the river hydropower generation $g$ $(g \in \mathcal{G}^{HRor})$ [MW]
\item [$P_{g}^{Acc}$] accumulated capacity of technology $p$ $(p \in \mathcal{P}^{*})$ [MW]
\item [$Q_{s}^{AccSE}$] accumulated storage capacity of electricity store $s$ ($s \in \mathcal{S}^{E}$) [MWh]
\item [$P_{l}^{AccL}$] accumulated capacity of line $l$ ($l \in \mathcal{L}$) [MW]
\item [$C^{R}$] operational cost of producing 1 kg hydrogen from SMRCCS [€/kg]
\item [$P^{HSea/HRor}_{gt}$] production profile of seasonal hydropower/run of the river hydropower in location $z$ period $t$ $(z \in \mathcal{Z}, t \in \mathcal{T})$ [MW]
\item [$V^{DHy}_{zt}$] hydrogen demand in region $z$ period $t$ $(z \in \mathcal{Z}, t \in \mathcal{T})$ [MW]
\item [$E^{R}$] emission factor of SMRCCS
\item [$V^{Acc}$] storage level, injection and withdrawal capacities of hydrogen storage and capacity of SMRCCS [kg]

\item [\textbf{Operational model variables}]
\item [$p_{gt}^{G}$] power generation of thermal generator $g$ in period $t$ ($g \in \mathcal{G}, t \in \mathcal{T}$) [MW]
\item [$p_{gt}^{ResG}$] power reserved of thermal generator $g$ for spinning reserve requirement in period $t$ ($g \in \mathcal{G}, t \in \mathcal{T}$) [MW]
\item [$p_{st}^{SE+}/p_{st}^{SE-}$] charge/discharge power of electricity store $s$ in period $t$ ($s \in \mathcal{S}^{E}, t \in \mathcal{T}$) [MW]
\item [$p_{st}^{ResSE}$] power reserved in electricity store $s$ for spinning reserve requirement in period $t$  ($s \in \mathcal{S}^{E}, t \in \mathcal{T}$) [MW]
\item [$q_{st}^{SE}$] energy level of electricity store $s$ at the start of period $t$ ($s \in \mathcal{S}^{E}, t \in \mathcal{T}$) [MWh]
\item [$p_{zt}^{GShed,l}$] generation shed for power $(l=P)$ and heat $(l=H)$ in location $z$ in period $t$ ($z \in \mathcal{Z}, t \in \mathcal{T}$) [MW]
\item [$p_{zt}^{Shed,l}$] load shed for power $(l=P)$ and heat $(l=H)$ in location $z$ in period $t$ ($z \in \mathcal{Z}, t \in \mathcal{T}$) [MW]
\item [$v_{zt}^{GShedHy}$] hydrogen production shed in location $z$ in period $t$ ($z \in \mathcal{Z}, t \in \mathcal{T}$) [kg]
\item [$v_{zt}^{ShedHy}$] hydrogen load shed in location $z$ in period $t$ ($z \in \mathcal{Z}, t \in \mathcal{T}$) [kg]
\item [$p_{lt}^{L}$] power flow in line $l$ in period $t$ ($l \in \mathcal{L}, t \in \mathcal{T}$) [MW]
\item [$p_{bt}^{BE}$] power consumption of electric boiler $b$ in period $t$ ($b \in \mathcal{B}^{E}, t \in \mathcal{T}$) [MW]
\item [$p^{F}_{ft}$] power generation of fuel cell $f$ in period $t$ ($f \in \mathcal{F}, t \in \mathcal{T}$) [MW]
\item [$p^{E}_{et}$] power consumption of electrolyser $e$ in period $t$ ($e \in \mathcal{E}, t \in \mathcal{T}$) [MW]
\item [$v^{SHy+}_{st}/v^{SHy-}_{st}$] injection/withdraw of hydrogen to (from) hydrogen storage $s$ in period $t$ ($s \in \mathcal{S}^{Hy}, t \in \mathcal{T}$) [kg]
\item [$v^{SHy}_{st}$] storage level of hydrogen storage $s$ in period $t$ $(s \in \mathcal{S}^{Hy}, t \in \mathcal{T})$ [kg]
\item [$v^{R}_{rt}$] hydrogen production of SMRCCS $r$ in period $t$ $(r \in \mathcal{R}, t \in \mathcal{T})$ [kg]
\item [$v^{LHy}_{lt}$] hydrogen flow in pipeline $l$ period $t$ $(l \in \mathcal{L}^{Hy}, t \in \mathcal{T})$
\item [$v_{vt}$] hydrogen injection, withdraw, storage level of hydrogen storage, and hydrogen production of SMRCCS in period $t$ $(v \in \mathcal{S}^{Hy} \cup \mathcal{R})$ [kg]
\end{description}

\end{document}